\documentclass[final,leqno]{siamltex}
\usepackage{amsmath,amssymb}
\usepackage{booktabs,dcolumn}
\usepackage{footnote}
\usepackage{float}
\usepackage[pdftex]{graphicx,color}
\usepackage{epstopdf}
\usepackage{makecell}
\usepackage{multicol}
\usepackage{enumerate}
\usepackage{float}
\usepackage{hyperref} \hypersetup{
   colorlinks=true,
 linktocpage=true,
 bookmarksnumbered=true
}
\usepackage{algorithm,algorithmic}

\newtheorem{thm}{Theorem}[section]
\newtheorem{rem}[thm]{Remark}

\newcommand{\name}{\text{SSI-TL}\,}
\newcommand{\nameexp}{\text{explicit-TL}\,}

\newcommand{\Li}{L^{\text{int}}}
\newcommand{\Ls}{L_c^{\text{exp}}}
\newcommand{\Lc}{L_c^{\text{imp}}}

\newcommand{\ie}{\emph{i.e.}}

\newcommand{\ud}{\mathrm{d}}
\newcommand{\prob}[1]{\mathrm{P}\left(#1\right)}
\newcommand{\expt}[1]{\mathrm{E}\left[#1\right]}

\newcommand{\var}[1]{\mathrm{Var}\left[#1\right]}

\newcommand{\norm}[1]{\left\|#1\right\|}
\newcommand{\abs}[1]{\left|#1\right|}
\newcommand{\rset}{\mathbb{R}}

\newcommand{\zset}{\mathbb{Z}}

\newcommand{\Ordo}[1]{{\mathcal{O}}\left(#1\right)}

\newcommand{\ordo}[1]{{o}\left(#1\right)}

\newcommand{\latt}{\mbox{$\zset_+^d$}}

\newcommand{\PERIOD}{.}
\newcommand{\COMMA}{,}

\newcommand{\LP}{\left(}
\newcommand{\RP}{\right)}

\newcommand{\algorithmfootnote}[2][\footnotesize]{  \let\old@algocf@finish\@algocf@finish  \def\@algocf@finish{\old@algocf@finish    \leavevmode\rlap{\begin{minipage}{\linewidth}
    #1#2
    \end{minipage}}  }}

\title{Multilevel Hybrid Split-Step Implicit  Tau-Leap}

\author{Chiheb Ben Hammouda\thanks{Computer, Electrical and Mathematical Science and Engineering,
 King Abdullah University of Science and Technology (KAUST),
 Thuwal, Saudi Arabia ({\tt chiheb.benhammouda@kaust.edu.sa}).}
        \and Alvaro Moraes\thanks{Computer, Electrical and Mathematical Science and Engineering,
 King Abdullah University of Science and Technology (KAUST),
 Thuwal, Saudi Arabia ({\tt alvaro.moraesgutierrez@kaust.edu.sa}).}
\and  Raul Tempone\thanks{Mathematical and Computer Science and Engineering Division,
 King Abdullah University of Science and Technology (KAUST),
 Thuwal, Saudi Arabia ({\tt raul.tempone@kaust.edu.sa}).}}

\makeatletter\begin{document}

\maketitle

\begin{abstract}
In biochemically reactive systems with small copy numbers of one or more reactant molecules, the dynamics is dominated by stochastic effects. To approximate those systems, discrete state-space and stochastic simulation approaches have been shown to be more relevant than continuous state-space and deterministic ones. 
In systems characterized by having simultaneously fast and slow timescales, existing discrete space-state stochastic path simulation methods,  such as the stochastic simulation algorithm (SSA) and the explicit tau-leap (\nameexp) method, can be very slow. 
Implicit approximations have been developed to improve numerical stability and provide  efficient simulation algorithms for those systems. Here, we propose an efficient Multilevel Monte Carlo (MLMC) method  in the spirit of the work by Anderson and Higham (2012). This method uses split-step implicit tau-leap (\name) at levels where the \nameexp method is not applicable due to numerical stability issues.
We present numerical examples that illustrate the performance of the proposed method.
 \end{abstract}



\pagestyle{myheadings}
\thispagestyle{plain}

\setcounter{tocdepth}{1}

\section{Introduction}\label{Introduction}

In his work, we extend the split-step backward-Euler method introduced in \cite{higham2002strong} to numerically solve stochastic differential equations (SDEs) driven by Brownian motion to the setting of SDEs driven by Poisson random measures \cite{cinlar2011probability,li2007analysis}.

We  focus on a particular class of continuous-time Markov chains named stochastic reaction networks (SRNs) (see Section \ref{sec:pjp} for a short introduction). SRNs are employed to describe the time evolution of biochemical reactions, epidemic processes \cite{anderson2015stochastic,SIR}, and transcription and translation in genomics and virus kinetics \cite{Genomics_SRN,srivastava2002stochastic},  among other important applications. 
Let $\mathbf{X}$ be an SRN taking values in $\latt$ and defined in the time-interval $[0,T]$, where $T>0$ is a user-selected final time. We aim to provide accurate estimations of  the expected value, 
$\expt{g(\mathbf{X}(T))}$,  where $g:\rset^d\to\rset$ is a given {smooth}  scalar observable of $\mathbf{X}$.

Many methods have been developed to simulate exact sample paths of SRNs;
for instance, the {stochastic simulation algorithm} (SSA) was
introduced by Gillespie in \cite{gillespie_ssa} and the {modified next reaction method} (MNRM) was proposed by Anderson in \cite{anderson:214107}. 
It is well known that pathwise exact realizations of SRNs may be  computationally very costly when some reaction channels have high reaction rates.
To overcome this issue, Gillespie \cite{gillespie_tau_leap} and Aparicio and Solari  \cite{Aparicio_tau_leap} independently proposed the explicit tau-leap (TL) method  (see Section \ref{sec:exp_tau})  to simulate approximate paths of  $\mathbf{X}$ by evolving the process with fixed time steps, keeping the reaction rates fixed within each time step. 
 
In this work, we address the problem of producing approximate-path simulations of SRNs for cases in which the set of reaction channels can be clearly classified into two subsets: fast and slow reaction channels (in the literature, this type of problem is called a \emph{stiff problem} \cite{abdulle2010chebyshev,Gillespie03}). 
In this context, the \nameexp method is hardly adequate because the time step required to maintain the numerical stability can be very small \cite{Cao_stab,rathinam} and, as a consequence, simulating  \nameexp paths of process $\mathbf{X}$ may become computationally expensive. For this purpose, different implicit schemes have been suggested, in particular, the {drift-implicit tau-leap} \cite{Gillespie03}, which uses an implicit discretization for the drift, together with an explicit discretization of the Poissonian noise. Our split-step implicit tau-leap (\name), method naturally produces values in the lattice, \latt. Hence, there is no need for an additional rounding step as in \cite{Gillespie03} (see Section \ref{sec:imp_tau}).

Other simulation schemes have been proposed to deal with situations with well-separated fast and slow timescales. 
Cao and Petzold \cite{Cao2005b} proposed a trapezoidal implicit tau-leap formula that is similar to the  trapezoidal rule for solving ordinary differential equations (ODEs). 
The trapezoidal implicit tau-leap scheme is numerically stable and does not possess the damping effect present in the drift-implicit tau-leap (Section 3 in \cite{Cao_stab} gives details on the damping effect). 
In \cite{Cao2005b}, the authors focus on the stability properties of a trapezoidal implicit tau-leap scheme.
They provide a posteriori comparisons of moments of order one and two with the SSA, the  \nameexp and the drift-implicit tau-leap methods. To the best of our  knowledge, no multilevel version of the trapezoidal tau-leap has been proposed.
In \cite{Ahn_implicit}, Ahn et al. developed fully implicit tau-leap schemes that implicitly treat both the drift and the noise of the Poisson variables. The construction of fully implicit tau-leap schemes is based on adapting weakly convergent discretizations of stochastic differential equations \cite{num_sol_sdes_Platen} to stochastic chemical kinetic systems. 
Constructing a multilevel version of the fully implicit tau-leap is still under investigation since it may involve a more complicated way of coupling two consecutive paths than does the \name presented here. In the same context, Rathinam et al. \cite{rathinam2007stiff} developed an explicit method that addresses stiffness in small number systems. The reversible-equivalent-monomolecular tau method (REMM-$\tau$) consists of  decomposing the system into  reversible pairs of reactions and  approximating bimolecular reversible pairs by suitable monomolecular reversible pairs. Updating the state of the system  involves the use of both binomial and Poisson random variables. The authors in \cite{rathinam2007stiff} showed that  REMM-$\tau$ performs better than the drift-implicit \cite{Gillespie03} and the trapezoidal implicit \cite{Cao_stab,Cao2005b} tau-leap methods  for small number stiff problems, since it avoids the need for an additional rounding step required in the latter methods. Howerver, REMM-$\tau$  does not have a natural multilevel Monte Carlo (MLMC) version (see Section \ref{Coupling Partial Implicit Tau-Leap Paths}). 
 Another type of explicit method that addresses stiffness was developed in \cite{abdulle2010chebyshev} by Abdulle et al. They proposed $\tau$-ROCK (Orthogonal Runge-Kutta Chebyshev) as an extension of  the multi-stage S-ROCK methods \cite{abdulle2007stabilized,S-ROCK_Abdulle,abdulle2008} for discrete stochastic processes and they derived schemes for processes with discrete Poisson noise. 
This method can have comparable stability properties to the implicit method while remaining explicit. 
The stability properties of $\tau$-ROCK can be controlled by increasing the tuning parameter, which is the stage number. 
Although the $\tau$-ROCK method has the advantage of avoiding the cost of solving nonlinear problems involved in \name, another cost is included that is related to the selection procedure for the optimal number of stages to obtain a desired stability domain. 
In \cite{abdulle2010chebyshev}, Abdulle et al. did not provide an explicit comparison between $\tau$-ROCK and drift-implicit tau-leap in terms of computational time and error control. 
Furthermore, to the best of our  knowledge, no multilevel version of $\tau$-ROCK has been developed for SRNs.
The multilevel estimator proposed here is based on the split-step implicit tau-leap (\name) rather than on $\tau$-ROCK for two reasons: the first  is that coupling two consecutive paths is simpler in the case of the \name than in the case of the one based on  $\tau$-ROCK, which is currently under investigation; the second is that simulating single paths with \name is more computationally efficient than with $\tau$-ROCK, especially in the case of large time steps (see Section \ref{sec:examples}).  Finally, we mention the so-called \textit{ integral tau methods} proposed in \cite{Rathinam_integral_tau}, in which the authors developed two variants of tau-leap methods that are adapted to stiff problems, produce nonnegative and lattice valued solutions,  and converge to the implicit Euler solution in the fluid limit. This is achieved in two steps: i) a deterministic implicit step that computes the drift update, and ii) a stochastic step based on the search for a set of feasible reaction counts that adds a linear combination of independent binomial and Poisson random variables. 
Nevertheless, those methods cannot be naturally extended to the MLMC paradigm.

We mention a few works on error analysis and convergence in implicit tau-leap schemes \cite{rathinam,rathinam2016convergence}. In \cite{rathinam}, Rathinam et al. showed that both the explicit and the implicit tau-leap methods are first-order convergent in all moments under the assumptions of bounded  state space and  linear propensity functions. This result was generalized for nonlinear propensity functions in \cite{li2007analysis} but only specifically for the  \nameexp method. This type of method has $1/2$ order  of  convergence in a strong sense. Discussions of the two existing scalings for the convergence analysis of tau-leap methods can be found in \cite{anderson2011error,hu2011weak}. More general results for the weak convergence were established in \cite{rathinam2016convergence} for  unbounded systems that possess certain moment growth bounds. Rathinam \cite{rathinam} showed that the main convergence result applies to any  tau-leap method provided that it produces integer-valued states and  satisfies  some condition related to  the moments, particularly implicit methods.

Inspired by the work of Anderson and Higham \cite{Anderson2012,GilesMLMC}, we extend the \name idea to the multilevel setting by introducing a multilevel hybrid \name   estimator (see Section \ref{sec:Multilevel Drift-Implicit Tau-leap}) with the aim of reducing the computational work needed to compute an estimate of $\expt{g(\mathbf{X}(T))}$ within a fixed tolerance, $TOL$,  with a given level of confidence, and for a class of systems characterized by having simultaneously fast and slow time scales. Our MLMC strategy couples two \name paths at the coarser levels of discretization until a certain interface level, $\Li$, where we start coupling paths using the  \nameexp method as indicated in \cite{Anderson2012}. In that sense, our strategy can be considered to be a hybrid algorithm. This strategy is specially relevant when $TOL$ is small, implying that the finest level of the multilevel \name, $L^{\text{imp}}$, is in the stability regime of the multilevel  \nameexp. For large values of $TOL$, our MLMC estimator reduces to a  multilevel \name  estimator. 

As noted above, this extension is relevant to systems with slow and fast timescales (see Section \ref{sec:examples} for numerical examples). In such situations, the multilevel estimator given in \cite{Anderson2012} is not computationally efficient due to the numerical stability constraint. 
In Section \ref{sec:examples}, we show that the multilevel  \name  method has the same order of computational complexity  as the multilevel  \nameexp, which is of $\Ordo{TOL^{-2} \log(TOL)^{2}}$  \cite{Anderson_Complexity}, but with a smaller multiplicative constant. We note that here our goal is to provide an estimate of $\expt{g(\mathbf{X}(T))}$ in the probability sense and not in the mean-square sense as in \cite{Anderson2012}.	

{
Recently, adaptive multilevel estimators were proposed to improve the performance of non-adaptive estimators \cite{Anderson2012} to simulate SRNs with markedly different timescales. 
In \cite{tau_control_variate}, Moraes et al. presented  an adaptive MLMC method that uses a low-cost, reaction-splitting heuristic to adaptively classify  the set of reaction channels into two subsets, fast and slow, in terms of level of activity. 
Based on their adaptive reaction-splitting technique and using a hierarchy of non-nested time meshes, they simulate the increments associated with high activity channels (fast reactions) using the  \nameexp method while those associated with low activity channels (slow reactions) are simulated using an exact method. 
Lester et al. \cite{lester2015adaptive} proposed an adaptive MLMC method strongly influenced by \cite{Anderson2012} that uses a time-stepping strategy based on the concept of \emph{the leap condition} introduced in \cite{Cao2006}. 
The idea of adaptivity provides possibilities for future research directions, for instance developing efficient, adaptive, hybrid multilevel estimators that would allow us to switch between using \name,  \nameexp and exact SSA within the course of a single sample path.}

The outline of the remainder of this work is as follows:
in Section \ref{sec:singlelevel}, we introduce the mathematical model of SRNs and give the necessary elements to simulate single-level \name paths in the context of SRNs.
In Section \ref{sec:Multilevel Drift-Implicit Tau-leap}, we introduce our multilevel hybrid \name  estimator. First, we review the MLMC method and then we show how to couple \name paths and also a \name path with another path simulated by the \nameexp method. This coupling procedure is the main building block for constructing the new MLMC estimator. This section ends with a presentation of the implementation details of the proposed estimator, which involves the procedure to select the levels and number of simulations per level and the switching rule from coupled \name to  \nameexp. In Section \ref{sec:examples}, we present some numerical experiments illustrating the performance of the proposed method. Finally, in Section \ref{conclusions}, we offer conclusions and suggest directions for future work.

\section{Single-Level \name Path-Simulation of Stochastic Reaction Networks}
\label{sec:singlelevel}
In this section, we briefly introduce the definition of stochastic reaction networks (SRNs). Then, we review the  \nameexp method. Finally, we present the \name scheme, which can be seen as an extension of the split-step backward-Euler method \cite{higham2002strong} to the context of SRNs.

\subsection{Stochastic Reaction Networks}
\label{sec:pjp}
We are interested in the time evolution of a homogeneously mixed  chemical reacting system described by the Markovian pure jump process,  $\mathbf{X}:[0,T]\times \Omega \to\zset_+^d$, where ($\Omega$, $\mathcal{F}$, $P$) is a probability space. In this framework, we assume that $d$ different species interact through $J$ reaction channels. 
The $i$-th component,  $X_{i}(t)$, describes the abundance of the $i$-th species present in the chemical system at time $t$. Hereafter, $\mathbf{v}^{T}$ denotes the transpose of the vector $\mathbf{v}$. The aim of this work is to study the time evolution of the state vector, 
\begin{equation*}
 \mathbf{X}(t) = (X_1(t), \ldots, X_d(t))^{T} \in
  \zset_+^d \PERIOD
\end{equation*}
Each reaction channel, $\mathcal{R}_j$, is a pair $(a_j, \boldsymbol{\nu}_{j})$, defined by its propensity function, $a_{j}:\rset^{d} \rightarrow \rset_{+}$, and its state change vector, $ \boldsymbol{\nu}_{j}=( \nu_{j}^{1},\nu_{j}^{2},..., \nu_{j}^{d})^{T}$, satisfying
\begin{align}\label{reaction_channel}
 \prob{\mathbf{X}(t+ \Delta t)=\mathbf{x}+ \boldsymbol{\nu}_{j} \bigm| \mathbf{X}(t)=\mathbf{x}}=a_{j}(\mathbf{x})\Delta t + \ordo{\Delta t}, \: j=1,2,...,J  \PERIOD 
\end{align}
Formula   \eqref{reaction_channel} states that the probability of observing a jump in the process, $\mathbf{X}$, from  state $\mathbf{x}$ to  state $\mathbf{x} +  \boldsymbol{\nu}_{j}$, a consequence of the firing of the j-th reaction, $\mathcal{R}_{j}$, during a small time interval, $(t, t + \Delta  t]$, is proportional to the length of the time interval, $\Delta  t$, with $a_{j}(\mathbf{x})$  as the constant of proportionality. 

We set $a_j(\mathbf{x}){=}0$ for those $\mathbf{x}$ such that $\mathbf{x}{+}\boldsymbol{\nu}_j\notin\latt$ (\emph{the non-negativity assumption}: the system can never produce negative population values). 

As a consequence of relation (\ref{reaction_channel}), process $\mathbf{X}$  is a continuous-time, discrete-space Markov chain that can be characterized by the random time change representation of Kurtz \cite{kurtzmp}: 
\begin{equation}
  \label{eq:exact_process}
\mathbf{X}(t)= \mathbf{x}_{0}+\sum_{j=1}^{J} Y_j \LP \int_0^t  a_{j}(\mathbf{X}(s)) \, \ud s \RP \boldsymbol{\nu}_j  \COMMA
\end{equation}
where $Y_j:\rset_+{\times} \Omega \to\zset_+$ are independent unit-rate Poisson
processes. 
Conditions on the set of reaction channels  can be imposed to ensure uniqueness  \cite{anderson2015stochastic} and to avoid explosions in finite time \cite{rathinam2013moment,gupta2014scalable,Engblo2014stab}

\subsection{The Explicit-TL Approximation}
\label{sec:exp_tau}
In this section, we  briefly review the  \nameexp approximation of the process, $\mathbf{X}$.  

The tau-leap is a pathwise-approximate method independently introduced in \cite{gillespie_tau_leap} and \cite{Aparicio_tau_leap} to overcome the computational drawback of exact methods, \ie, when many reactions fire during a short time interval. It can be derived from the random time change representation of Kurtz  \eqref{eq:exact_process} by approximating the integral $\int_{t_i}^{t_{i+1}} a_{j}(\mathbf{X}(s)) \ud s $ by $a_j(\mathbf{X}(t_i))\,(t_{i+1}-t_i)$, \ie, using the forward-Euler method with a time mesh  $\{t_{0}=0, t_{1},...,t_{N}= T\}$. In this way, the  \nameexp approximation of $\mathbf{X}$ should satisfy for $k\in\{1,2,\ldots,N\}$: 
\begin{equation*}\label{approx}
\mathbf{Z}^{exp} (t_{k}) = \mathbf{x}_{0}+\sum_{j=1}^{J} Y_{j} \LP  \sum_{i=0}^{k-1} a_{j}(\mathbf{Z}^{exp}(t_i))(t_{i+1}-t_{i}) \RP   \boldsymbol{\nu}_{j} \PERIOD
\end{equation*}
By considering a uniform time mesh of size $\tau$, we can simulate a path of $\mathbf{Z}^{exp}$ as follows. Let $\mathbf{Z}^{exp} (t_0) := \mathbf{x}_{0}$ and define 
\begin{equation}\label{explicit}
\mathbf{Z}^{exp} (t_k):=\mathbf{z}+\sum_{j=1}^{J} \mathcal{P}_{j}(a_{j}(\mathbf{z}) \tau)  \boldsymbol{\nu}_{j} \COMMA
\end{equation}
iteratively, where $\mathbf{z}=\mathbf{Z}^{exp} (t_{k-1})$ and $\mathcal{P}_{j}(r_j) $ are independent Poisson random variables with respective rates, $r_j$. Notice that the  \nameexp path, $\mathbf{Z}^{exp}$, is defined only at the points of the time mesh, but it can be naturally extended to $[0,T]$  as a piecewise constant path by defining $\mathbf{Z}^{exp} (t_{k-1}+h) := \mathbf{Z}^{exp} (t_{k-1}), \forall\,\, 0<h<\tau$.

\subsubsection{Numerical Stability of the Explicit-TL Method}
The numerical stability of the  \nameexp method is treated in \cite{rathinam} for the case of linear propensities, \ie,  $\: a_j(\mathbf{X}) = \mathbf{c}^{T}_j \mathbf{X}$, where $\mathbf{c}_{j} \in \mathbb{R}^{d}$. 
In this particular case, taking expectations  of relation \eqref{explicit} conditional on $\mathbf{z}$  results in
\begin{align}\label{explicit_cond_expextation}
\expt{ \mathbf{Z}^{exp} (t_k)\mid \mathbf{z} }= (1+\tau A) \mathbf{z} \COMMA
\end{align}
where the $d \times d$ matrix $\mathbf{A}$ is given by $ \mathbf{A}=\sum\limits_{j=1}^{J} \boldsymbol{\nu}_{j} \mathbf{c}_{j}^{T} $.

Taking the expectation of \eqref{explicit_cond_expextation}, we obtain
\begin{align}\label{explicit_cond_expextation_expextation}
\expt{ \mathbf{Z}^{exp} (t_k)}= (1+\tau A) \expt{\mathbf{z}} \COMMA
\end{align}
 which yields 
 \begin{align}
\expt{ \mathbf{Z}^{exp} (t+N \tau)}= (1+\tau A)^{N} \expt{\mathbf{Z}^{exp} (t)}.
\end{align}

Then, from \eqref{explicit_cond_expextation_expextation}, the expectation of the  \nameexp method is asymptotically stable if $\tau$ satisfies
\begin{align*}
\mid 1+\tau \lambda_{i}(\mathbf{A}) \mid <1, \: \: i=1,\dots, d,
\end{align*}
where $\{ \lambda_{i} \}_{i=1}^{d}$ are the eigenvalues of the matrix $\mathbf{A}$. 

In the general case, the numerical stability limit of the  \nameexp method can be computed by a linearized stability analysis of the forward-Euler method applied to the corresponding deterministic ODE model \cite{Gillespie03}. In the case of systems having simultaneously fast and slow timescales, this stability limit can be very small, implying an expensive computational cost for simulation. To overcome this issue, the drift-implicit tau-leap idea \cite{Gillespie03} has been proposed. We should note that using implicit methods is only relevant in the case of negative eigenvalues. In fact,  having a positive  eigenvalue implies a rapid
change in the solution,  which requires the use of a small time step.  Using explicit methods is more appropriate under such conditions.

\subsection{The \name Approximation}
\label{sec:imp_tau}
In this section, we define  $\mathbf{Z}^{imp}$, the \name  approximation of the process, $\mathbf{X}$. 

The  \nameexp scheme (\ref{explicit}), where $\mathbf{z}=\mathbf{Z}^{exp} (t)$,  can be rewritten as follows:
\begin{align}\label{explicit2}
\mathbf{Z}^{exp} (t+\tau) &=\mathbf{z}+\sum_{j=1}^{J} \mathcal{P}_{j} \LP a_{j}(\mathbf{z}) \tau \RP \boldsymbol{\nu}_{j} \nonumber\\
 &=\mathbf{z}+\sum_{j=1}^{J} \LP \mathcal{P}_{j} \LP a_{j}(\mathbf{z}) \tau \RP - a_{j}(\mathbf{z}) \tau + a_{j}(\mathbf{z}) \tau \RP  \boldsymbol{\nu}_{j} \nonumber\\
 &=  \mathbf{z}+\sum_{j=1}^{J} a_{j}(\mathbf{z})  \tau  \boldsymbol{\nu}_{j} +\sum_{j=1}^{J} \LP \mathcal{P}_{j} \LP a_{j}(\mathbf{z}) \tau \RP - a_{j}(\mathbf{z}) \tau \RP  \boldsymbol{\nu}_{j} \PERIOD
\end{align}

Let us denote the second and third quantities on the right-hand side of (\ref{explicit2}) by the drift and the zero-mean noise, respectively. The idea of the \name  method is to take only the drift part as implicit while the noise part is left explicit. Let us define
 $\mathbf{z}=\mathbf{Z}^{imp} (t)$ and define 
$\mathbf{Z}^{imp} (t+\tau)$ through the following  two steps:
\begin{align}\label{partial_implicit}
\mathbf{y} &= \mathbf{z}+\sum_{j=1}^{J}  a_{j} \LP \mathbf{y}\RP    \tau \boldsymbol{\nu}_{j}\text{ (Drift-Implicit step)}\\\nonumber
\mathbf{Z}^{imp} (t+\tau) &=  \mathbf{y} + \sum_{j=1}^{J} \LP \mathcal{P}_{j}(a_{j}(\mathbf{y}) \tau)-  a_{j}(\mathbf{y})\tau\RP   \boldsymbol{\nu}_{j}\\
\nonumber
&=  \mathbf{z} + \sum_{j=1}^{J}  \mathcal{P}_{j}(a_{j}(\mathbf{y}) \tau)   \boldsymbol{\nu}_{j}\text{ (Tau-leap step)} \PERIOD
\end{align}

\begin{algorithm}[h!]
\caption{Algorithm of the \name method}.   

\label{alg:imp1}
\begin{algorithmic}[1]
	\STATE Fix $ \tau>0$. Set $\mathbf{Z}(0)=\mathbf{x}_{0}$, $t_{0}=0$, $n=0$ and repeat the following until $t_{n} \geq T$
	\STATE Set $t_{n+1}=t_{n}+\tau$. If $t_{n+1}\geq T$, set $t_{n+1}=T$ and $\tau=T-t_n$
		\STATE Implicit step: $ \mathbf{Y}=\mathbf{Z}(t_n)  +\sum\limits_{j=1}^{J}   a_{j}(\mathbf{Y})  \tau \boldsymbol{\nu}_{j} $
		\STATE Explicit step: $\mathbf{Z}(t_{n+1})=\mathbf{Z}(t_n)+ \sum\limits_{j=1}^{J}  \left( \text{Poisson}(a_{j}(\mathbf{Y}) \tau)\right)  \boldsymbol{\nu}_{j}$
    \STATE Set $n \leftarrow n+1$ 
\end{algorithmic}
\end{algorithm}
Algorithm \ref{alg:imp1} implements \eqref{partial_implicit}. 

Notice that the \name method naturally produces values of $Z^\text{imp}$ in the lattice, \latt, without the need for a final rounding step as in the drift-implicit method \cite{Gillespie03}, which goes as follows: first compute an intermediate state according to relation \eqref{intermediate_state}:
\begin{align}\label{intermediate_state}
 \mathbf{Z}^{'}&=  \mathbf{z}+\sum_{j=1}^{J} a_{j}( \mathbf{Z}^{'}) \:  \tau  \: \boldsymbol{\nu}_{j} +\sum_{j=1}^{J} \LP \mathcal{P}_{j} \LP a_{j}(\mathbf{z}) \tau \RP - a_{j}(\mathbf{z}) \tau \RP  \boldsymbol{\nu}_{j} \PERIOD
\end{align}
Then  approximate the actual number of firings, $K_j(\mathbf{z}, \tau)$, of reaction channel $R_j$ in the time interval $(t, t + \tau]$ by the integer-valued random variable, $\widetilde{K}_j(\mathbf{z}, \tau)$, defined by
\begin{align}\label{firing_rounding}
\widetilde{K}_j(\mathbf{z}, \tau)&=  \left[a_{j}( \mathbf{Z}^{'}) \tau + \mathcal{P}_{j} \LP a_{j}(\mathbf{z}) \tau \RP - a_{j}(\mathbf{z}) \tau \right]  \COMMA
\end{align}
where $[x]$ denotes the nearest nonnegative integer to $x$.  Finally, update the state at time $t+\tau$ such that 
$\mathbf{Z} (t+\tau) =  \mathbf{z} + \sum_{j=1}^{J} \widetilde{K}_j(\mathbf{z}, \tau)   \boldsymbol{\nu}_{j}.$ This rounding procedure in relation \eqref{firing_rounding} introduces an error when simulating the path.

\section{Multilevel Hybrid \name}
\label{sec:Multilevel Drift-Implicit Tau-leap}
The goal of this section is to define our multilevel hybrid \name estimator. 
First, we quickly review the MLMC method as proposed by Giles \cite{GilesMLMC}. 
Then, we show how to couple two \name paths associated with two nested time meshes. 
This coupling is the main building block for constructing our novel multilevel hybrid \name estimator. 
Finally, we present the switching rule that we use to select the time mesh in which we move from the \name method to the  \nameexp one.
For simplicity, in this section, we consider only uniform time meshes in the interval $[0,T]$.

\subsection{The Multilevel Monte Carlo (MLMC) Method}
\label{sec:MLMC)}
Let $\mathbf{X}$ be a stochastic process and $g: \rset ^{d} \rightarrow \rset$ a smooth scalar observable. 
Let us assume that we want to approximate $\expt{g(\mathbf{X}(T))}$, but instead of sampling directly from $\mathbf{X}(T)$, we sample from $\mathbf{Z}_{h}(T)$, which are random variables generated by an approximate method with step size $h$.  
Let us assume also that the variates $\mathbf{Z}_{h}(T)$ are generated with an algorithm with weak order, $\Ordo{h}$, \ie, $\expt{g(\mathbf{X}(T))- g(\mathbf{Z}_{h}(T) )}= \Ordo{h}$.
 
Let $\mu_{N}$ be the standard Monte Carlo estimator  of $\expt{g(\mathbf{Z}_{h}(T))}$ defined by
\begin{equation}\label{MC estimator}
\mu_{N} :=\frac{1}{N}\sum_{n=1}^{N} f(\mathbf{Z}_{h,[n]}(T))\COMMA
\end{equation}
where $\mathbf{Z}_{h,[n]}(T)$ are independent and distributed as $\mathbf{Z}_{h}(T)$. 

Consider now the following decomposition of the global error:
\begin{equation}
\expt{g(\mathbf{X}(T))} - \mu_{N}= \LP \expt{g(\mathbf{X}(T))- g(\mathbf{Z}_{h}(T) )}  \right)+ \left( \expt{g(\mathbf{Z}_{h}(T) )} - \mu_{N} \RP \PERIOD
\end{equation}
To have the desired accuracy, $TOL$, it is sufficient to take  $h=\Ordo{TOL}$ so that the first term on the right is $\Ordo{TOL}$ and, by the Central Limit Theorem,  impose $N=\Ordo{TOL^{-2}}$ so that the statistical error given by the second term on the right is $\Ordo{TOL}$ \cite{duffie1995efficient}.  
As a consequence, the expected total computational work is $\Ordo{TOL^{-3}}$.

The MLMC estimator, introduced by Giles \cite{GilesMLMC} allows us to reduce the total computational work from $\Ordo{TOL^{-3}}$ to $\Ordo{TOL^{-2} \log(TOL)^{2}}$. 
The basic idea of MLMC is to generate,  and couple in an intelligent manner, paths with different step sizes, which results in
\begin{itemize}
\item[i)] Stochastically coordinated sequences of paths having different step sizes, where the paths with large step sizes are computationally less expensive than those with very small step sizes.
\item[ii)] A small variance of the difference between two coupled paths with fine step sizes, implying significantly fewer samples in the estimation.
\end{itemize}
We can construct the MLMC estimator as follows: consider a hierarchy of nested meshes of the time interval $[0,T]$, indexed by $\ell=0, 1,\dots, L$. We denote by $h_{0}$ the step size used at level $\ell=0$. The size of the subsequent time steps for levels $\ell \geq 1$ are given by $ h_{\ell}=M^{-\ell} h_{0}$, where $M{>}1$ is a given integer constant. In this work, we take $M = 2$. To simplify the notation,  hereafter $\mathbf{Z}_{\ell}$ denotes the approximate process generated using a step size of $h_{\ell}$.  

Consider now the following telescoping decomposition of $\expt{g(\mathbf{Z}_{L}(T))}$:
\begin{equation}\label{MLMC1}
\expt{g(\mathbf{Z}_{L}(T))}= \expt{g(\mathbf{Z}_{0}(T))}+ \sum_{\ell=1}^{L} \expt{g(\mathbf{Z}_{\ell}(T))- g(\mathbf{Z}_{\ell-1}(T))} \PERIOD
\end{equation} 
By defining 
\begin{align}\label{estimators}
\begin{cases} 
\hat{Q}_{0}:= \frac{1}{N_{0}} \sum\limits_{n_{0}=1}^{N_{0}} g(\mathbf{Z}_{0,[n_{0}]}(T)) \\ 
\hat{Q}_{\ell}:= \frac{1}{N_{\ell}} \sum\limits_{n_{\ell}=1}^{N_{\ell}}  \LP g(\mathbf{Z}_{\ell,[n_{\ell}]}(T))-g(\mathbf{Z}_{\ell-1,[n_{\ell}]}(T)) \RP \COMMA \\
\end{cases}
\end{align}
we arrive at the unbiased MLMC estimator, $\hat{Q}$, of  $\expt{g(\mathbf{Z}_{L}(T))}$:
\begin{equation}\label{MLMCest}
\hat{Q}:= \sum\limits_{\ell=0}^{L} \hat{Q}_{\ell}\PERIOD
\end{equation}

 We note that the key point here is that both $\mathbf{Z}_{\ell,[n_{\ell}]}(T)$ and $\mathbf{Z}_{\ell-1,[n_{\ell}]}(T)$ are sampled using different time discretizations but with the same generated randomness. If we simulate the paths with a method having a strong error of order $1/2$, \ie, $\expt{ \mid \mid \mathbf{Z}_{h}(T)-\mathbf{X}(T)\mid \mid^{2}}= \Ordo{h}$, and we assume that $g$ is  Lipschitz, \ie, $\exists \: c \geq 0 \: \: \text{s.t} \:  \mid g(\mathbf{X})- g(\mathbf{Y})\mid \leq c \mid \mid \mathbf{X} -\mathbf{Y} \mid \mid $, it is straightforward to see that 
 $$ \var{g(\mathbf{Z}_{\ell}(T))- g(\mathbf{Z}_{\ell-1}(T))} = \Ordo{h_{l}}.$$
Therefore, by setting $N_{\ell}=\Ordo{ TOL^{-2} \times L \times h_{\ell}}$, we obtain $ \var{\hat{Q}}= \Ordo{TOL^{2}}$ but with a total computational complexity of $\Ordo{TOL^{-2} \log(TOL)^{2}} $, which makes the MLMC estimator better than the standard MC estimator \eqref{MC estimator} for computational purposes.

\subsection{Coupling Two \name Paths}\label{Coupling Partial Implicit Tau-Leap Paths}

Algorithm \ref{alg:MLMC_imp}, inspired by \cite{Anderson2012}, shows how to couple two \nameexp paths. Note that $d$ is the number of species and $J$ is the number of reactions.

\begin{algorithm}[h!]
\caption{Simulation of two coupled \name paths.}
\label{alg:MLMC_imp}
\begin{algorithmic}[1]
	\STATE Fix $ h_\ell > 0$  and set $h_{\ell-1} = 2 \times h_\ell$. Set $\mathbf{Z}_{\ell}(0) = \mathbf{Z}_{\ell-1}(0) = \mathbf{x}_{0}$, $t_{0} = 0$, $n = 0$. 
		\WHILE{$t_{n} < T$}
	\STATE Implicit step for the finer level: 
	
	  Solve $ \widetilde{\mathbf{Z}}_{\ell}(t_{n}) = \mathbf{Z}_{\ell}(t_{n})+\sum\limits_{j=1}^{J}    a_{j}(\widetilde{\mathbf{Z}}_{\ell}(t_{n})) h_{\ell} \boldsymbol{\nu}_{j}$
	\STATE Implicit step for the coarser level:
	   \IF{$(n \: mod \: 2) 	= 0 $}  		
				\STATE  Solve $   \widetilde{\mathbf{Z}}_{\ell-1}(t_{n}) = \mathbf{Z}_{\ell-1}(t_{n})+\sum\limits_{j=1}^{J}    a_{j}( \widetilde{\mathbf{Z}}_{\ell-1}(t_{n})) h_{\ell-1} \boldsymbol{\nu}_{j}$
	   \ENDIF
	\FOR {$j{=}1$ \textbf{to} $J$}
	\STATE $A_{3(j-1)+1}= a_{j}(\widetilde{\mathbf{Z}}(t_{n})) \wedge a_{j}(\widetilde{\mathbf{Z}}_{\ell-1}(t_{n}))$
	\STATE $A_{3(j-1)+2}= a_{j}(\widetilde{\mathbf{Z}}(t_{n}))-A_{3(j-1)+1} $
	\STATE $A_{3(j-1)+3}=a_{j}(\widetilde{\mathbf{Z}}_{\ell-1}(t_{n}))-A_{3(j-1)+1} $
	\STATE $\Lambda_{3(j-1)+1}= \: \text{Poisson} \: (A_{3(j-1)+1} h_\ell) $
	\STATE $\Lambda_{3(j-1)+2}= \: \text{Poisson} \: (A_{3(j-1)+2} h_\ell)$
	\STATE $\Lambda_{3(j-1)+3}= \: \text{Poisson} \: (A_{3(j-1)+3} h_\ell) $
	\ENDFOR
	\STATE State updating (explicit step)
	\begin{itemize}
	\item[i)] Set $\boldsymbol{\Gamma}_{\ell}= \boldsymbol{\nu} \otimes [1 \: 1 \: 0]$ and $\boldsymbol{\Gamma}_{\ell-1}= \boldsymbol{\nu} \otimes [1 \: 0 \: 1]$ 
	($A \otimes B$ refers to the Kronecker product of the matrices $A$ and $B$, therefore 
	$\boldsymbol{\Gamma}_{\ell}$ and $\boldsymbol{\Gamma}_{\ell-1}$ are $d \times 3 J$ matrices. Observe that $A$ and $\Lambda$ are of size $3J\times 1$.	

	\item[ii)] Update $\mathbf{Z}_{\ell}(t_{n+1}) = \mathbf{Z}_{\ell}(t_{n})+ h_{\ell} \boldsymbol{\Gamma}_{\ell} \Lambda$
	\item[iii)] Update $\mathbf{Z}_{\ell-1}(t_{n+1}) = \mathbf{Z}_{\ell-1}(t_{n})+ h_{\ell} \boldsymbol{\Gamma}_{\ell-1} \Lambda$
	\end{itemize}
	\STATE $t_{n+1}=t_{n}+h_{\ell}$
	\STATE $n \leftarrow n+1$,
	\ENDWHILE
\end{algorithmic}
\end{algorithm}

We define our multilevel hybrid \name estimator in the next section based on Algorithm \ref{alg:MLMC_imp}.

\begin{rem}
For solving the implicit steps in Algorithm \ref{alg:MLMC_imp}, we use the classic Newton method.
As an initial guess, we select $\mathbf{X}(0)$. Regarding the number of iterations, we ideally would like to compute an approximate solution, $\mathbf{y^*}$, of $\mathbf{y}$ such that the distribution of $\text{Poisson}(a_j(\mathbf{y^*})\,h)$ is close to the distribution of $\text{Poisson}(a_j(\mathbf{y})\,h)$ for all $j \in \{1,2,\ldots,J\}$. 
In practice, we usually check that the difference between two consecutive iterations, 
$h\,\max_j \norm{\nu_j}_{\infty}\,\abs{a_j(\mathbf{y_{k-1}^*})- a_j(\mathbf{y_{k}^*})}$, is below a certain tolerance.
\end{rem}

\begin{rem}
It is well known that tau-leap methods can produce negative
population numbers. Methods for avoiding negative population numbers can
roughly be divided into three classes: i) the pre-leap check technique \cite{Gillespie_avoiding_negative_steps,Cao2006,alvaro_chernoff}; ii) the post-leap check technique \cite{anderson_postleap}; and iii) the modification of the Poisson distributed increments with bounded increments from binomial or multinomial distributions \cite{Burrage_Binomial_leap}. The pre-leap check technique  computes the largest possible time step satisfying some leap criterion. It is often based on controlling the relative change in the propensity function before taking the step.  
The post-leap check procedure is applied when  a step results in a negative  population value. In this case, it retakes a shorter step, conditioned on already sampled data from the failed step, to avoid sampling bias. 
Finally, the third technique consists of  projecting  negative values of species numbers to zero at each step  and continues the path simulation until $T$. In the future, we intend to investigate these techniques in the context of the \name method.
\end{rem}    
    
\subsection{Definition of the Multilevel Hybrid  \name Estimator}
\label{sec:Hybrid Unbiased Tau-leap MLMC Estimator}
Our hybrid MLMC estimator uses the  multilevel \name method only at the coarser levels and then, starting from a certain interface level, $\Li$, it switches to the  multilevel  \nameexp method as given in \cite{Anderson2012}.  This strategy is specially relevant when $TOL$ is small, implying that the finest level of the multilevel \name, $L^{\text{imp}}$, is in the stability regime of the multilevel \nameexp. 
On the contrary, for sufficiently large values of the tolerance, $TOL$, our estimator is reduced to a multilevel \name estimator. 
For simplicity, let us consider a family of uniform time meshes of $[0,T]$, with size $h_{\ell} = 2^{- \ell} T$. 
Let $\Lc$  and $\Ls$ be the coarsest levels in which the \name and the explicit methods are respectively stable. 
In the class of problems we are interested in, we have the relation $h_{\Ls} \ll h_{\Lc}$,  which means that  $\Lc$ is much coarser than $\Ls$.

Rewriting \eqref{MLMCest} in our context, we define our multilevel hybrid  \name estimator as
\begin{equation}\label{MLMCest_hybrid}
\hat{Q}:= \hat{Q}_{\Lc}+\sum\limits_{\ell=\Lc+1}^{\Li-1} \hat{Q}_{\ell}+ \hat{Q}_{\Li} +\sum\limits_{\ell=\Li+1}^{L} \hat{Q}_{\ell}\COMMA
\end{equation}

where 
\begin{align}\label{estimators_hybrid}
\begin{cases} 
\hat{Q}_{\Lc}:= \frac{1}{N_{i,\Lc}} \sum\limits_{n=1}^{N_{i,\Lc}} g(\mathbf{Z}^{imp}_{\Lc,[n]}(T)) \\ 
\hat{Q}_{\ell}:= \frac{1}{N_{ii,\ell}} \sum\limits_{n_{\ell}=1}^{N_{ii,\ell}}  \LP g(\mathbf{Z}^{imp}_{\ell,[n_{\ell}]}(T))-g(\mathbf{Z}^{imp}_{\ell-1,[n_{\ell}]}(T)) \RP, \: \:  \Lc+1 \leq \ell \leq  \Li-1  \\
\hat{Q}_{\Li}:= \frac{1}{N_{ie,\Li}} \sum\limits_{n=1}^{N_{ie,\Li}}  \LP g(\mathbf{Z}^{exp}_{\Li,[n]}(T))-g(\mathbf{Z}^{imp}_{\Li-1,[n]}(T)) \RP \\
\hat{Q}_{\ell}:= \frac{1}{N_{ee,\ell}} \sum\limits_{n_{\ell}=1}^{N_{ee,\ell}}  \LP g(\mathbf{Z}^{exp}_{\ell,[n_{\ell}]}(T))-g(\mathbf{Z}^{exp}_{\ell-1,[n_{\ell}]}(T)) \RP  , \: \:  \Li+1 \leq \ell \leq  L  \PERIOD \\
\end{cases}
\end{align}
{
\begin{rem}\label{hyb_est_coupling}
According to \eqref{estimators_hybrid}, computing our MLMC tau-leap estimator \eqref{MLMCest_hybrid} may require three types of coupling: i) coupling two consecutive \name paths (see section \ref{Coupling Partial Implicit Tau-Leap Paths}), ii) coupling  a \name path with an  \nameexp path; this coupling is made in a similar way as in Algorithm \ref{alg:MLMC_imp} but instead of an implicit step for the finer level, we do an explicit step, and iii) coupling two consecutive  \nameexp paths \cite{Anderson2012}.
\end{rem}}

From relations \eqref{MLMCest_hybrid} and \eqref{estimators_hybrid}, we notice that our MLMC estimator requires the specification of the coarsest discretization level, $\Lc$,  the interface level, $\Li$, the finest level of discretization, $L$, and the number of samples per level defined by $\mathbf{N} := \{N_{i,\Lc}, \{N_{ii,\ell}\}_{\ell=\Lc+1}^{\Li-1},N_{ie,\Li},\{ N_{ee,\ell}\}_{\ell=\Li+1}^{L}\}$.

\subsubsection{On the Selection of the Coarsest Discretization Level}
We start by defining the criterion that we use to choose the value of the coarsest discretization level, $\Lc$. In fact, to ensure the numerical stability of our MLMC estimator, two conditions must be satisfied: the first one ensures the stability of a single path, which is related to the coarsest level of discretization, $L_{c}$, and  which can be determined by a linearized stability analysis of the backward-Euler method applied to the deterministic ODE model corresponding to our SRN system \cite{Gillespie03}. The second one ensures the stability of the variance of the coupled paths of our MLMC estimator and can be expressed by $\var{g(\mathbf{Z}_{\Lc+1}){-}g(\mathbf{Z}_{\Lc})} \ll \var{g(\mathbf{Z}_{\Lc})}$.

\subsubsection{On the Selection of the Finest Discretization Level}
The total number of levels, $L$, and the set of the number of samples per level, $\mathbf{N}$,  are selected to satisfy the accuracy constraint, 
$\prob{\abs{\expt{g(\mathbf{X}(T))} - \hat Q}<TOL}>1-\alpha$  (typically we choose $\alpha = 0.05$),  with near-optimal expected computational work. Here, $TOL$ is a user-selected tolerance.

{
The total error can be split into bias and statistical error such that
\begin{align}
\abs{\expt{g(\mathbf{X}(T))} - \hat Q} \leq   \underbrace{\abs{\expt{g(\mathbf{X}(T)) - \hat Q}}}_{\text{Bias}}+  \underbrace{\abs{\expt{\hat Q} - \hat Q}}_{\text{Statistical error}} \PERIOD
\end{align}

If we use a splitting parameter, $\theta \in (0,1)$, satisfying
\begin{align*}
TOL=  \underbrace{(1-\theta) \,\,TOL}_{\text{Bias tolerance}}+   \underbrace{\theta \,\,TOL}_{\text{Statistical error tolerance}} \COMMA
\end{align*}
then, using the same idea introduced in \cite{collier_continuation_MLMC}, the MLMC algorithm should bound the bias and the statistical error as follows:
\begin{align}
\abs{\expt{g(\mathbf{X}(T)) - \hat Q}} & \leq (1-\theta)\,\,  TOL \COMMA  \label{MLMC_error_bound_bias} \\ 
\abs{\expt{\hat Q} - \hat Q} & \leq \theta\,\,  TOL \COMMA \label{MLMC_error_bound_stat}
\end{align}
where the latter bound should hold with probability $1{-}\alpha$. As stated in \cite{collier_continuation_MLMC}, relation \eqref{MLMC_error_bound_stat} can be achieved if we impose

\begin{align}\label{Stat_error_control}
\var{\hat Q} \leq \left( \frac{\theta \,\,TOL}{C_{\alpha}} \right)^{2} 
\end{align}
for some given confidence parameter, $C_{\alpha}$, such that $\Phi(C_{\alpha}) = 1 - {\alpha}/{2}$ ; here,  $\Phi$ is the cumulative distribution function of a standard normal random variable (see \cite{collier_continuation_MLMC} for details).

In our numerical examples, we choose $\theta = 0.5$. An optimal split between bias and statistical errors can be reached by assuming the dependence of $\theta$ on $TOL$ and applying the Continuation Multilevel Monte Carlo (CMLMC) algorithm, introduced in \cite{collier_continuation_MLMC}. Investigating the optimal split in this context will be left as  future work.}

In our problem, the finest discretization level, $L$, is determined by satisfying relation \eqref{MLMC_error_bound_bias} for $\theta = 0.5$, implying
\begin{equation}
\abs{\text{Bias}(L) :=   \expt{g(\mathbf{X}(T))- g(\mathbf{Z}_{L}(T))}} < \frac{TOL}{2}\PERIOD
\end{equation}
In our numerical experiments, we use the following approximation (see \cite{GilesMLMC}): 
\begin{equation}
\text{Bias}(L)\approx \expt{g(\mathbf{Z}_{L}(T)) -g(\mathbf{Z}_{L-1}(T))}\PERIOD
\end{equation}
Therefore, to determine the value of $L$, we need to have estimates of the bias for different levels of discretization, $\ell$.
In our numerical examples if $L\leq \Ls$: then we estimate the bias of the difference by coupling two \name paths. Otherwise,  the estimation is based on coupling two  \nameexp paths.

\begin{rem}
Due to the presence of the  large kurtosis (see Section 1 of \cite{alvaro_ML}) problem, as the level, $\ell$, increased,  we obtained $\mathbf{Z}_{\ell}(T) = \mathbf{Z}_{\ell-1}(T)$ in most of our simulations, while observing differences only in a very small proportion of the simulated coupled paths. For that reason, we extrapolate the bias and the variance of the consecutive differences obtained from the coarsest levels. Establishing dual-weighted formulae for the \name method, like the one proposed in \cite{alvaro_ML} for the multilevel \nameexp method, is  still under investigation.
\end{rem}


\subsubsection{On the Selection of the Number of Simulated Tau-Leap Paths and the Optimal Value of the Interface Level}
Given $L$,  we are now interested in determining a near-optimal number of samples per level given by $\mathbf{N}$ and also the optimal value of the interface level, $\Li$. For this purpose, we define $W_{\Li}$ to be the expected computational cost of the MLMC estimator given that the interface level is $\Li$, given by
\begin{align}\label{eq:work}
  W_{\Li} & := C_{i,\Lc}  N_{i,\Lc} h_{\Lc}^{-1}  + \sum_{\ell=\Lc+1}^{\Li-1} C_{ii,\ell} N_{ii,\ell} h_{\ell}^{-1} + C_{ie,\Li}   N_{ie,\Li} h_{\Li}^{-1} \nonumber \\
  &  +\sum_{\ell=\Li+1}^{L} C_{ee,\ell} N_{ee,\ell} h_{\ell}^{-1} \COMMA
 \end{align} 
 where $C_{i}$,  $C_{ii}$, $C_{ie}$ and $C_{ee}$ are, respectively, the  expected computational costs of simulating a single \name step,  a coupled \name step, a coupled \name / \nameexp step and  a coupled  \nameexp step.  
 These costs can be modeled as shown in Table \ref{Cost model.}.
\begin{table}[H]
\centering
\begin{tabular}{c | c }
cost &  cost model\\
\hline \hline
$C_{i}$ &  $C_{P}+C_{N}$ \\
$C_{ii}$ &  $\gamma (C_{P}+C_{N})$\\
$C_{ie}$ &  $\eta (C_{P}+C_{N})$ \\
$C_{ee}$ &  $\beta C_{P} $ \\
\hline
\end{tabular}
\caption{ Here $C_{P}$ represents the cost of generating $J$ Poisson random variables and $C_{N}$ represents the cost of the Newton iterations. The constants $\gamma$, $\eta$ and $\beta$ depend on the coupling and they are also machine-dependent quantities. We estimated these constants in our numerical experiments and we found that for Example 1: $\gamma \approx \beta \approx 2.6$ and $\eta \approx 1.8$ and for Example 2: $\gamma \approx \beta \approx 2.8$ and $\eta \approx 2$.  Examples are in the numerical section (Section \ref{sec:examples}). }
\label{Cost model.}
\end{table}

Now we fix the interface level, $\Li$, and see $W_{\Li}$ as a function of the number of samples, $\mathbf{N}$.

The first step is to solve
\begin{align}\label{opt_sampl_hyb1}
\begin{cases} 
\underset{\mathbf{N}}{\operatorname{min}} \: W_{\Li}(\mathbf{N})\\ 
s.t. \:   C_{\alpha} \sqrt{\sum\limits_{\ell=\Lc}^{L} N_{\ell}^{-1} V_{\ell}} \leq  \frac{TOL}{2} \COMMA
\end{cases}
\end{align}
where $ V_\ell = \var{g(\mathbf{Z}_{\ell}(T))-g(\mathbf{Z}_{\ell-1}(T)}$ is estimated, like the bias, by the extrapolation of the sample variances obtained from the coarsest levels. The confidence level  is usually taken as $C_{\alpha}=1.96$ (for $\alpha=0.05$ and assuming a Gaussian distribution of the estimator \cite{collier_continuation_MLMC}, see figures \ref{fig:normality_test} and \ref{fig:normality_test2} in Section \ref{sec:examples}). 

 The constraint in the optimization problem (\ref{opt_sampl_hyb1}) aims to control the statistical error of the MLMC estimator, given by relation \eqref{Stat_error_control} for $\theta=
0.5$, since the variance of this estimator, $\hat Q$, is given by $\sum\limits_{\ell=\Lc}^{L} N_{\ell}^{-1} V_{\ell}$ (the bias has been already controlled). 
 
\begin{rem}
In our numerical experiments, we notice that we obtain a large number of needed samples for the coarsest level due to the large variance at this discretization level. Reducing this variance, and hence the corresponding computational cost, can be achieved by using the control-variate technique. In our work, we used the idea proposed in \cite{tau_control_variate}, where the authors introduced a novel control-variate technique based on the stochastic time change representation by Kurtz \cite{kurtzmp}, which dramatically reduces the variance of the coarsest level of the multilevel Monte Carlo estimator  at a negligible computational cost.
\end{rem}


Let us denote $\mathbf{N^*}(\Li)$ as the solution for \eqref{opt_sampl_hyb1}.
Then, the optimal value of the switching parameter, $\Li$, should be chosen to 
minimize the expected computational work; that is, the value ${\Li}^*$ that solves
\begin{align*}
\begin{cases} 
\underset{\Li}{\operatorname{min}} \: W_{\Li}(\mathbf{N^*}(\Li)) \\ 
s.t. \: \Ls \leq \Li \leq L   \PERIOD 
\end{cases}
\end{align*}
 
By analyzing the cost per level, we can conclude that the lowest computational cost is most likely to be achieved for $\Li=\Ls $, \ie, the same level in which the  \nameexp is stable, see Figures \ref{fig:cost_hybrid_exp1_tol001} and \ref{fig:cost_hybrid_exp2_tol001} in Section \ref{sec:examples}.
 To motivate this selection of the interface level, we can write 
\begin{align}
W_{\Li+1}^{*} - W_{\Li}^{*} &=  C_{i}  \LP N_{i,\Lc}^{*}(\Li{+}1)-N_{i,\Lc}^{*}(\Li)\RP h_{\Lc}^{-1} 
\nonumber \\
&+ \sum_{\ell=\Lc+1}^{\Li-1} C_{ii,\ell} ( N_{ii,\ell}^{*}(\Li{+}1)-N_{ii,\ell}^{*}(\Li)) h_{\ell}^{-1} \nonumber \\
&+\sum_{\ell=\Li+2}^{L} C_{ee,\ell} ( N_{ee,\ell}^{*}(\Li{+}1)-N_{ee,\ell}^{*}(\Li))  h_{\ell}^{-1} \nonumber \\
&+  h_{\Li{+}1}^{-1} (  C_{ie} N_{ie,\Li{+}1}^{*}(\Li+1) - C_{ee} N_{ee,\Li{+}1}^{*}(\Li) )\nonumber \\
&+  h_{\Li}^{-1} ( C_{ii} N_{ii,\Li}^{*}(\Li{+}1)- C_{ie} N_{ie,\Li}^{*} (\Li) ) \PERIOD \nonumber 
\end{align}

In our numerical examples  in Section \ref{sec:examples}, we notice that $N_{\ell,ii}^{*}(\Li+1) \approx N_{\ell,ii}^{*}(\Li), \: \:  \: \Lc \leq\ell \leq \Li-1$ and $N_{\ell,ee}^{*}(\Li+1) \approx N_{\ell,ee}^{*}(\Li), \: \:  \Li+2 \leq\ell \leq L$,  which implies that 
\begin{align}\label{hybrid_work_analysis}
W_{\Li+1}^{*} - W_{\Li}^{*}  &\approx \underbrace{h_{\Li+1}^{-1} (C_{ie} N_{ie,\Li+1}^{*}(\Li+1) - C_{ee} N_{ee,\Li+1}^{*}(\Li) )}_{c_1}\nonumber\\
& - \underbrace{h_{\Li}^{-1} (  C_{ie} N_{ie,\Li}^{*} (\Li) -C_{ii} N_{ii,\Li}^{*}(\Li+1))}_{c_2} |\COMMA
\end{align}
see Figures \ref{fig:samples_hybrid_exp1_tol001} and \ref{fig:samples_hybrid_exp2_tol001}.

Our numerical experiments in Section \ref{sec:examples} show that $c_1>c_2$, which
means that we have $W_{\Li+1}> W_{\Li}$ and thus $\Li \approx \Ls$, see Figures  \ref{fig:cost_hybrid_exp1_tol001} and \ref{fig:cost_hybrid_exp2_tol001}.
{
\begin{rem}
Let us compare the expected computational work of the multilevel \nameexp and the multilevel  hybrid \name estimators, denoted as $W_{MLMC}^{\text{exp}} (TOL)$ and $W_{MLMC}^{\text{hyb}}(TOL)$, respectively. 
First, notice that the total number of \name levels, $\Li$, can be very large depending on how stiff our problem is, but it can not grow to infinity as $TOL \to 0$ because when the time mesh is sufficiently fine, we switch to the explicit method.
Second, for a time mesh of size $h_\ell$, we know that $V_\ell = \Ordo{h_\ell}$, where the constant depends on if the coupling is between \name or explicit paths as described in Table \ref{Cost model.}. As a consequence, for $TOL \to 0$, there is no advantage in using the \name method; that is:
\begin{align*}
\underset{TOL \rightarrow 0}{\lim} \frac{W_{MLMC}^{\text{hyb}} (TOL)}{W_{MLMC}^{\text{exp}}(TOL)}=1 \PERIOD
\end{align*} 
But let us observe that for strongly stiff problems and reasonable values of $TOL$, we do not expect to switch to the explicit method, \ie, only the first two terms of our estimator, $\hat Q$, are non-zero. 
In those cases,
 \begin{align*}
\frac{W_{MLMC}^{\text{hyb}} (TOL)}{W_{MLMC}^{\text{exp}}(TOL)} \ll1 \PERIOD
\end{align*} 
\end{rem}
}

{\begin{rem}
In  Section \ref{sec:examples}, we show that the asymptotical complexity of our multilevel \name estimator is  $ \Ordo{TOL^{-2} \log(TOL)^2}$ (see Figures \ref{fig:Comparison of the expected total work for the different method}  and \ref{fig:Comparison of the expected total work for the different method_exp2}). 
In principle, this complexity can be improved up to $\Ordo{ TOL^{-2}}$ (the optimal one in the context of Monte Carlo sampling). To achieve this optimal complexity,  it is enough to  couple  with a pathwise exact method at the bottom level, but, due to the high computational work of this alternative, the computational complexity of  $\Ordo{ TOL^{-2}}$ may not be observed in stiff problems.
\end{rem}}
 \section{Numerical Examples}
\label{sec:examples}
In this section, we present our numerical results from tests of the performance of our multilevel hybrid \name estimator on two examples.
\subsection{Example 1: The Decaying-Dimerizing Reaction Example}
The decaying-dimerizing reaction \cite{gillespie_tau_leap} consists of three species, $S_1$, $ S_2$, and $S_3$, and four reaction channels:
\begin{align}\label{dimerisation}
S_{1} &\overset {c_1}{\rightarrow} 0 \nonumber \\
S_{1}+ S_{1} &\overset {c_2}{\rightarrow} S_2  \nonumber \\
S_2 &\overset {c_3}{\rightarrow} S_{1}+ S_{1} \nonumber \\
S_{2} &\overset {c_4}{\rightarrow} S_3
\end{align}
We take the same values for the different parameters as in \cite{Gillespie03}: 
$c_1=1$, $c_2=10$, $c_3=10^{3}$ and $c_4=10^{-1}$. The stoichiometric vectors are $\boldsymbol{\nu}_{1}=(-1,0,0)^{T}$, $\boldsymbol{\nu}_{2}=(-2,1,0)^{T}$, $\boldsymbol{\nu}_{3}=(2,-1,0)^{T}$, and $\boldsymbol{\nu}_{4}=(0,-1,1)^{T}$. The corresponding propensity functions are
\begin{equation}
a_{1}(\mathbf{X}) = X_1,\: a_{2} (\mathbf{X})= 5 X_1(X_1 - 1),\: a_{3} (\mathbf{X})= 1000 X_2, \: a_{4}(\mathbf{X}) = 0.1 X_2,
\end{equation}
where $X_i$ denotes the number of molecules of species $S_i$ and the initial condition is $\mathbf{X}(0)=(400,798,0)^{T}$ [molecules]. We consider the final time, $T= 0.2$ seconds. In the following numerical experiments, we are interested in approximating $\expt{X_{3}(T)}$.  This setting implies that the stability limit of the \nameexp is $\tau_{\text{exp}}^{\text{lim}} \approx 2.3 \times 10^{-4}$ (computed using a linearized stability analysis of the forward-Euler method applied to the deterministic ODE model corresponding to this system \cite{Gillespie03}).

\begin{figure}[H]
 \begin{center}
\includegraphics[scale=0.35]{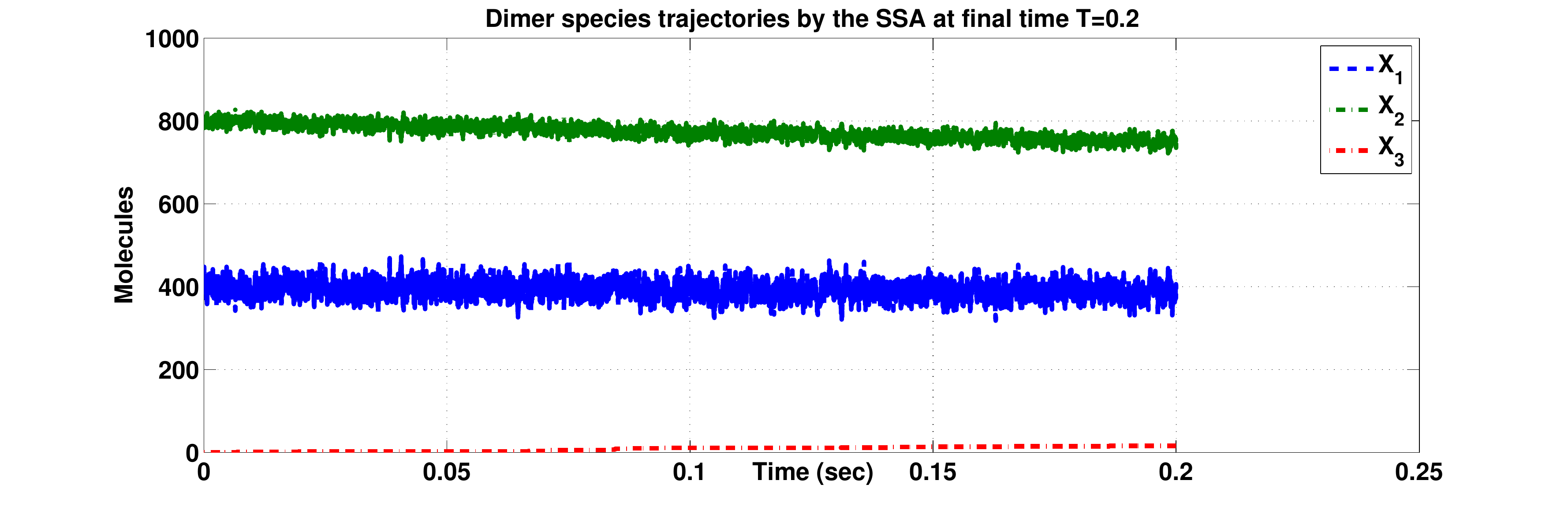}
\caption {Trajectories of dimer species trajectories simulated by SSA at the final time, $T=0.2$. This system has multiple timescales. 
The step size, $\tau^{\text{exp}}$, is therefore taken to be extremely small to ensure the numerical stability of the  \nameexp method ($\tau_{\text{exp}}^{\text{lim}} \approx 2.3 \times 10^{-4})$.}
\label{fig:Dimer_SSA}
\end{center}
\end{figure}

Figure (\ref{fig:Dimer_SSA}) shows the trajectories of dimer species with the reaction set (\ref{dimerisation}) solved with the original SSA. 

\subsubsection{Coupling Results}
\begin{figure}[H]
 \begin{center}
\includegraphics[scale=0.35]{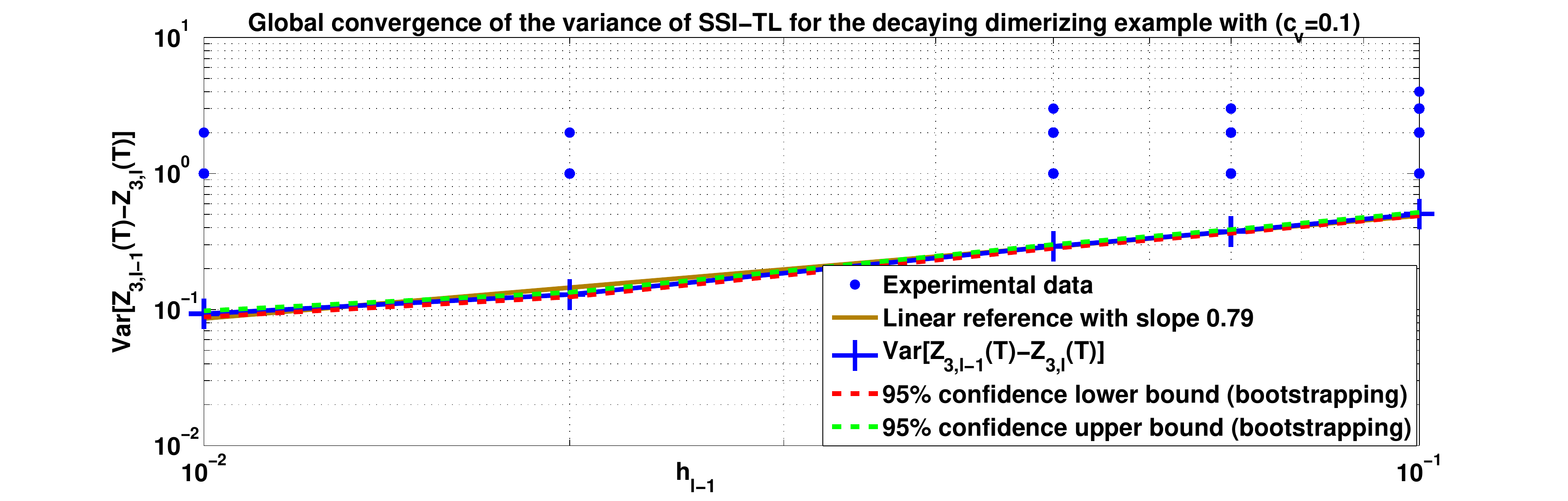}
\caption{Global Convergence of the variance of  \name for the decaying dimerizing example. The blue dots correspond to the observed data (each realization of $Z_{3,\ell}-Z_{3,\ell-1}$). We fit  $\var{Z_{3,\ell}-Z_{3,\ell-1}}$ using a bootstrapping technique \cite{Bootstrap} with a  coefficient of variation, $c_v=0.1$. There are some outliers due to the large kurtosis problem (see Section 1 of \cite{alvaro_ML}).}
\label{fig:Global_Strong_Convergence_implicit_decayin_dimerizing3}
  \end{center}
\end{figure}

\begin{figure}[H]
 \begin{center}
\includegraphics[scale=0.35]{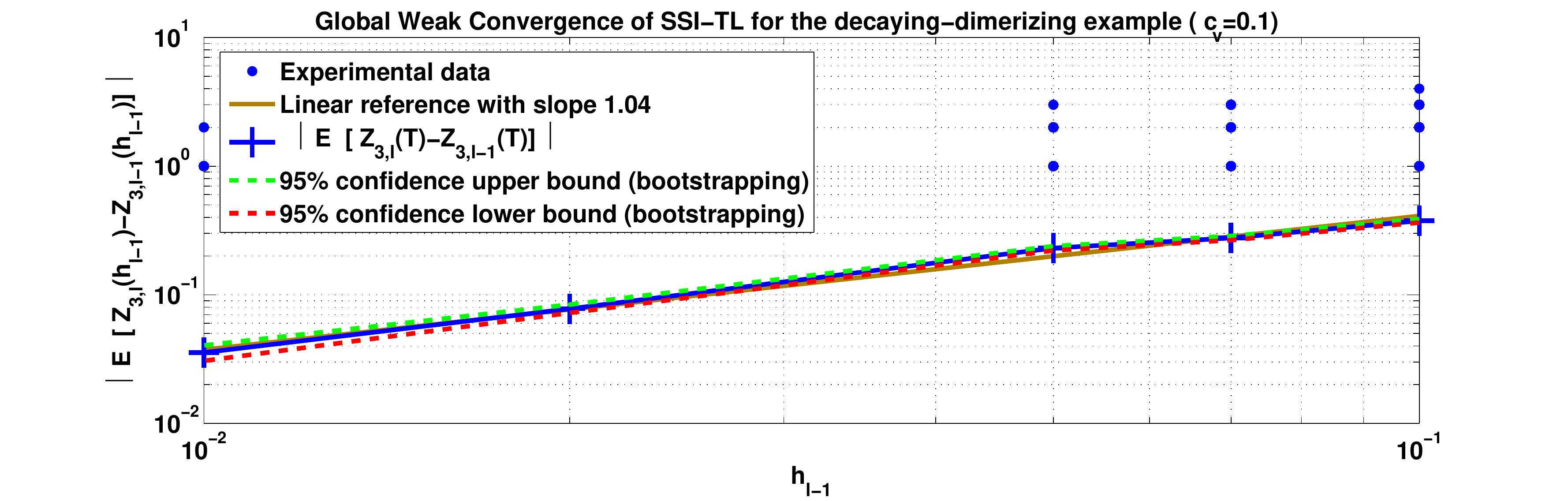}
\caption{Global weak convergence \name for the decaying-dimerizing example. Here, we fit $\expt{Z_{3,\ell}-Z_{3,\ell-1}}$ using a bootstrapping technique with a  coefficient of variation, $c_v=0.1$.}
\label{fig:Global_Weak_Convergence_implicit_decayin_dimerizing3}
  \end{center}
\end{figure}

The first step to construct our multilevel hybrid \name estimator is to check that the coupling procedure is correct in terms of convergence properties. The convergence tests that we did for the coupling algorithm indicate that the  global weak convergence and the global convergence of the variance are   approximately of   order $1$ (see Figures \ref{fig:Global_Strong_Convergence_implicit_decayin_dimerizing3} and \ref{fig:Global_Weak_Convergence_implicit_decayin_dimerizing3}). These results validate our coupling algorithm, which can be used to construct the corresponding MLMC estimator.

\subsubsection{Multilevel Hybrid SSI-TL Results ($L \leq \Ls$)}

Now, we compare the performance of our multilevel hybrid \name algorithm, reduced to the multilevel \name algorithm   when $L \leq \Ls$,  to the performance of the  multilevel \nameexp and   MC \name  methods. First, we should mention that to ensure the numerical stability of the  multilevel \nameexp estimator, we should be restricted to the stability limit given by $\tau_{\text{exp}}^{\text{lim}} \approx 2.3 \times 10^{-4}$. However, in the implicit case, we do not have any  stability restriction. Therefore, in our numerical experiments, we started with $\Ls=10$ for the  multilevel \nameexp estimator ($\Ls$ satisfies $h_{\Ls}= 2^{-\Ls} T < 2.3  \times 10^{-4}$). We also varied the prescribed tolerance, $TOL$, to investigate its effect on the performance of the tested methods. 

The finest level of discretization, $L$, given in Figure 6, is obtained by extrapolating the weak error and deduced from Figures 4 and 5 . The value of the finest level satisfies $Bias(L) < \frac{TOL}{2}$. From Figure \ref{fig:The finest level (L) function of tolerance _exp1}, we notice that  to achieve the same weak error, we use seven more levels for the \nameexp compared to those for \name. As a consequence, the weak error for the implicit scheme is a factor of $2^7 = 128$ smaller than for the explicit method. Even though this high gain in discretization is deteriorated by the additional cost of Newton iterations at each level, it is still one of the main factors that explains the great advantage of the multilevel \name over the explicit version as we will see later. 
\begin{figure}[H]
\centering
 \includegraphics[scale=0.35]{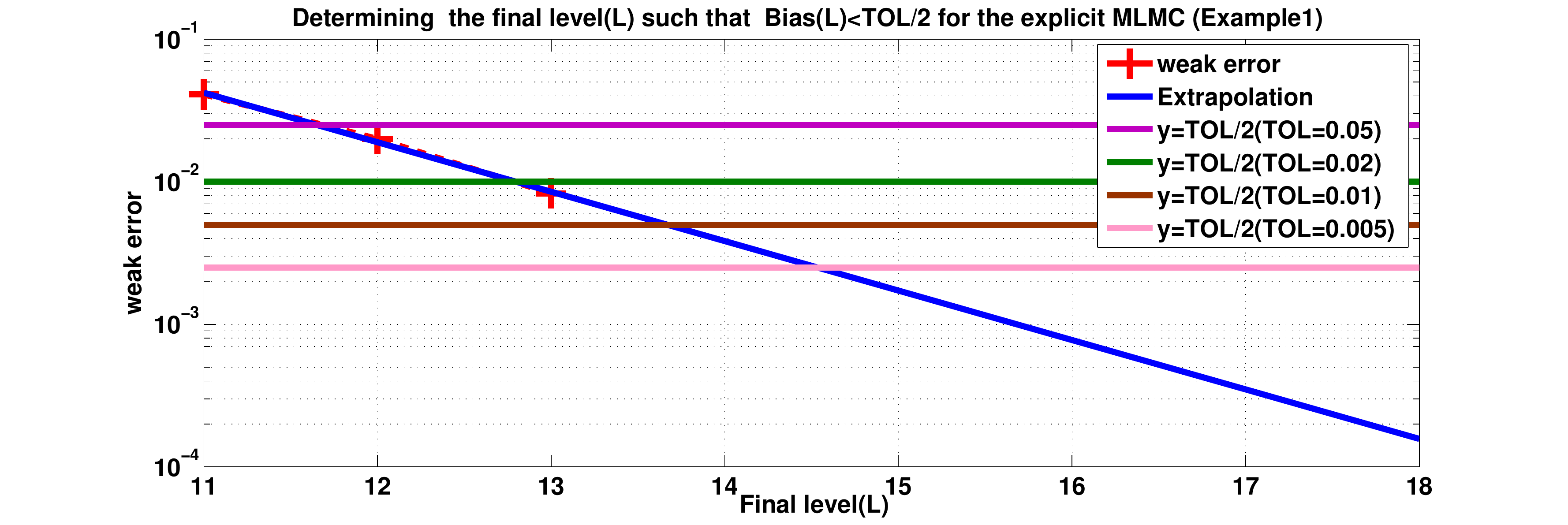}
  \caption [Determining the final level (L) such that  $Bias(L)<TOL/2$ for the  multilevel \nameexp (Example 1).]{Determining the final level (L) such that  $Bias(L)<TOL/2$ for the  multilevel \nameexp (Example 1). We estimate the bias for small values of level $\ell$ and then we extrapolate to obtain final level $L$.}
  \label{fig:1exp1}
\end{figure}
\begin{figure}[H]
\centering
 \includegraphics[scale=0.35]{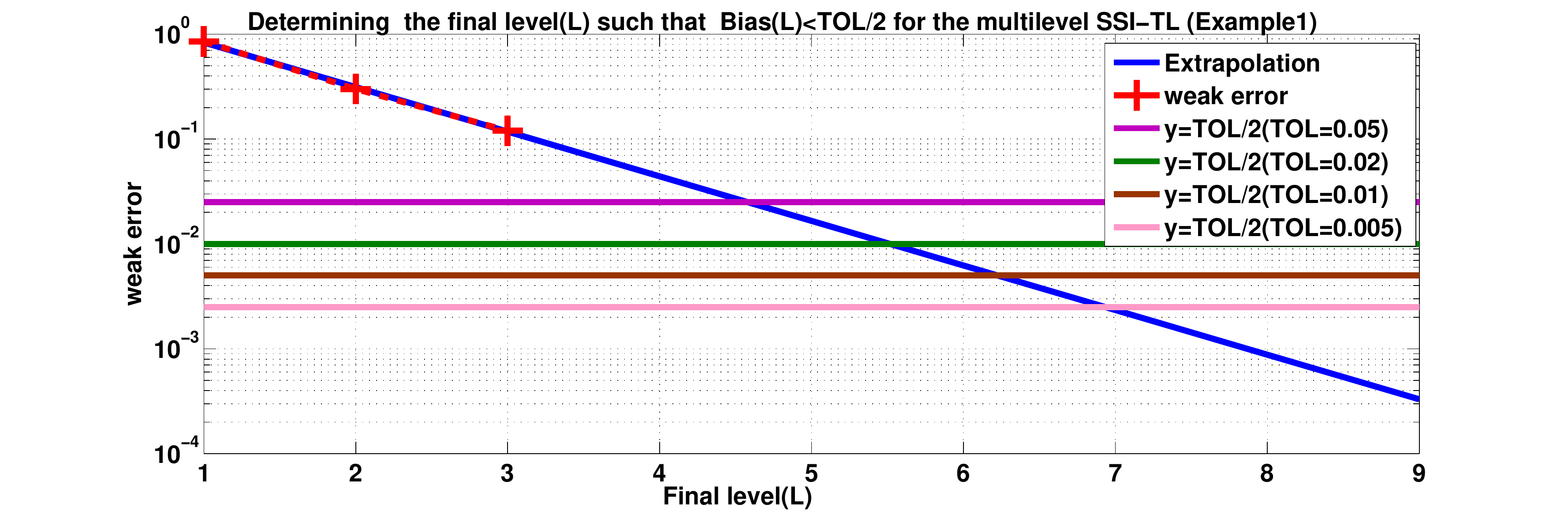}
  \caption{Determining the final level(L) such that  $Bias(L)<TOL/2$ for the multilevel \name (Example 1).}
  \label{fig:2exp1}
\end{figure}
\begin{figure}[H]
 \begin{center}
\includegraphics[scale=0.35]{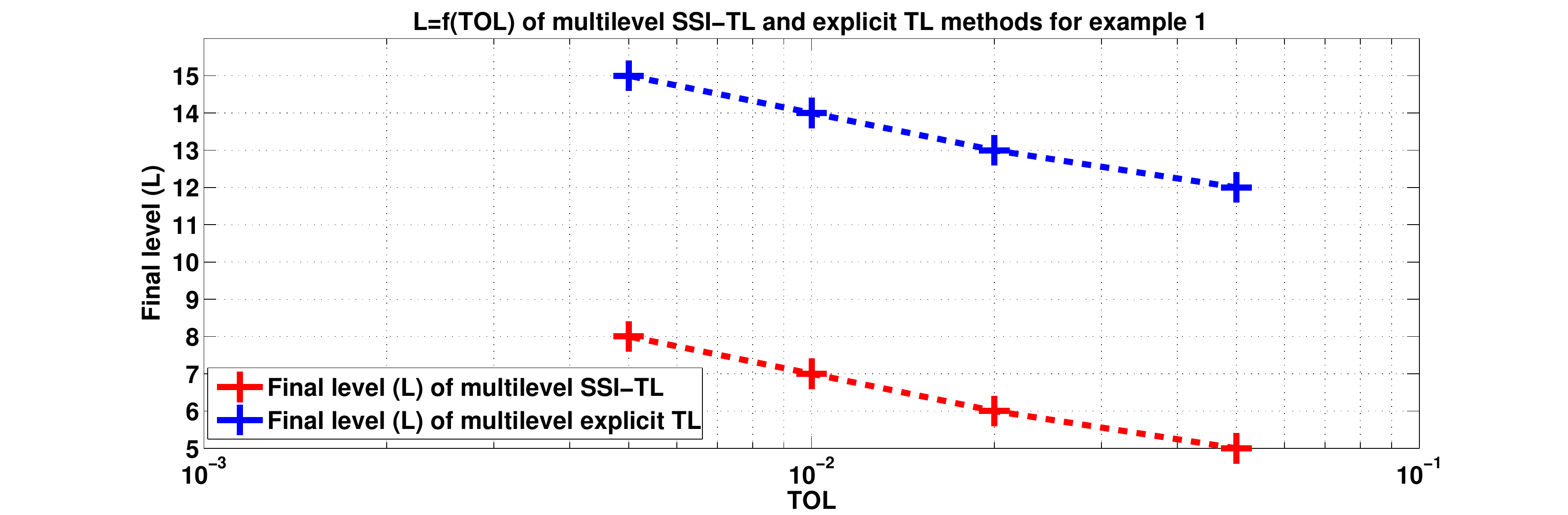}
\caption{The finest level (L) as a function of tolerance, $TOL$, for multilevel \name and multilevel \nameexp methods for Example 1. These values are obtained as a result of figures (\ref{fig:1exp1} and \ref{fig:2exp1}).}
\label{fig:The finest level (L) function of tolerance _exp1}
 \end{center}
\end{figure}

 The variances, $\{ V_{\ell} \}_{\ell=\Lc=0}^{L^{\text{imp}}}$, used in the optimization problem (\ref{opt_sampl_hyb1}) are obtained by the extrapolation of the estimated variances of the first three levels and given in  Table \ref{table:var_extrapol_1}.
\begin{table}[H]
\centering
\begin{tabular}{l | c | c | c }
 \multicolumn{2}{c}{\nameexp}   &   \multicolumn{2}{c}{\name}    \\
\hline \hline
level &     variance  &  level      &   variance    \\
\hline \hline
$11$  &    $4.1 \times 10^{-2}$ &  $1$   &  $3.1 $    \\
$12$ & $1.9 \times 10^{-2}$  & $2$ &   $1.7  $ \\
$13$ &    $8.2 \times 10^{-3}$   &  $3$ &   $9 \times 10^{-1}$   \\
$14$ &    $ 3.7 \times 10^{-3}$    &   $4$ &    $4.9 \times 10^{-1}$  \\
$15$ &  $1.6 \times 10^{-3}$ &   $5$ &    $2.6  \times 10^{-1}$ \\
$16$ & $7 \times 10^{-4}$  &    $6$ &   $1.4 \times 10^{-1}$  \\
$17$ & $3 \times 10^{-4}$  &     $7$ & $7.7 \times 10^{-2}$  \\
$18$ & $10^{-4}$ &   $8$ &   $  4.1 \times 10^{-2}$  \\
\hline
\end{tabular}
\caption{Extrapolation for the estimated variances, $ \var{Z_{\ell,3}-Z_{\ell-1,3}} \: (\ell=11,12,13)$ for the multilevel  \nameexp and $(\ell=1,2,3)$  for multilevel \name (Example 1).}
\label{table:var_extrapol_1}
\end{table}

\begin{figure}[H]
 \begin{center}
\includegraphics[scale=0.35]{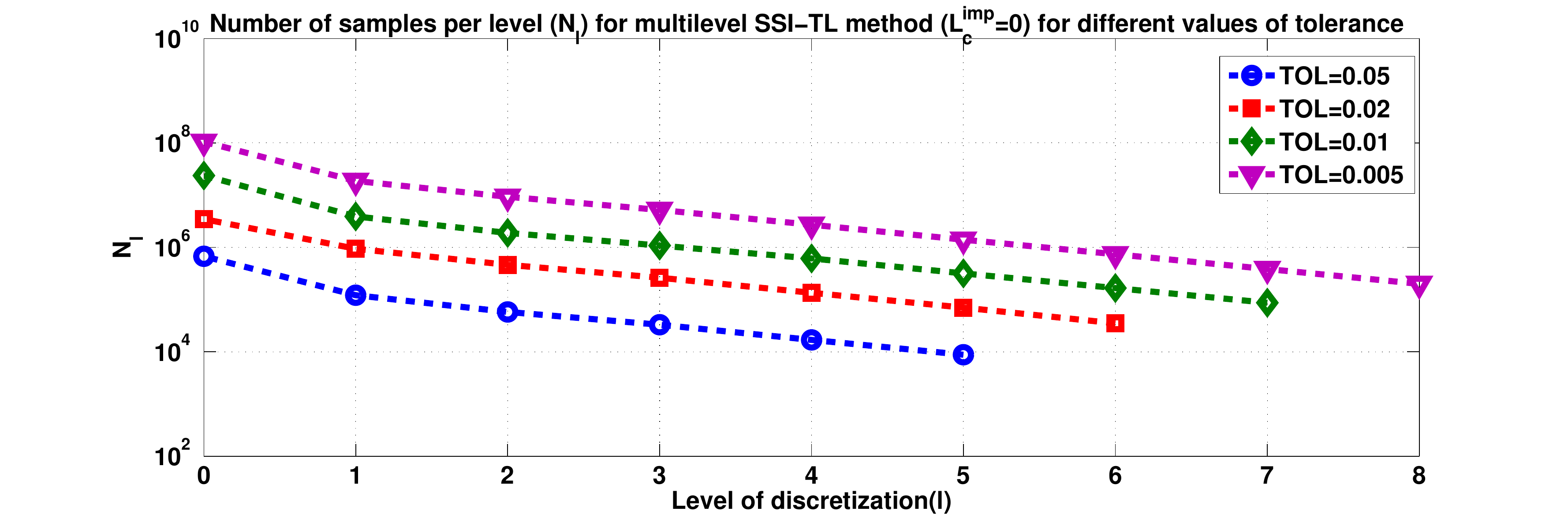}
\caption{Optimal number of samples per level for the multilevel \name method ($\Lc=0$) for Example 1. These values are obtained by solving the optimization problem (\ref{opt_sampl_hyb1}). The optimal number of samples increases as we decrease the tolerance. It decreases with respect to the level, $\ell$, due to the decrease of both $V_{\ell}$ and $h_{\ell}$. Note that the value of the finest level, $L$, the maximum level of time step refinement, increases as the value for the tolerance decreases.}
\label{fig:number_samples_implicit_1rst_example}
  \end{center}
\end{figure}

 To obtain the optimal number of samples per level for the multilevel \name estimator, we followed the procedure described in Section \ref{sec:Hybrid Unbiased Tau-leap MLMC Estimator} and our numerical results are presented in Tables \ref{Samples per level for the MLMC explicit tau-leap for the decaying-dimerizing example} and \ref{Samples per level for the drift implicit MLMC tau-leap for the decaying-dimerizing example} in the Appendix. These tables show the optimal number of samples for each estimator per level and for different values of tolerances, $TOL=\{0.05,\: 0.02, \:0.01,\:  0.005 \}$. From the previous tables and Figure \ref{fig:number_samples_implicit_1rst_example}, we see that the optimal number of samples increases as we decrease the tolerance. It decreases with respect to the level, $\ell$, due to the decrease in both $V_{\ell}$ and $h_{\ell}$. We also mention that the value of the finest level, $L$, the maximum level of time step refinement, increases as the value of the tolerance decreases. The control-variate technique that we used allowed us to reduce the number of samples at the coarsest level by a factor of six (see Tables \ref{Samples per level for the drift implicit MLMC tau-leap for the decaying-dimerizing example} and \ref{Samples per level for the drift implicit MLMC tau-leap for the decaying-dimerizing example control_variate}).

\begin{figure}[H]
 \begin{center}
\includegraphics[scale=0.35]{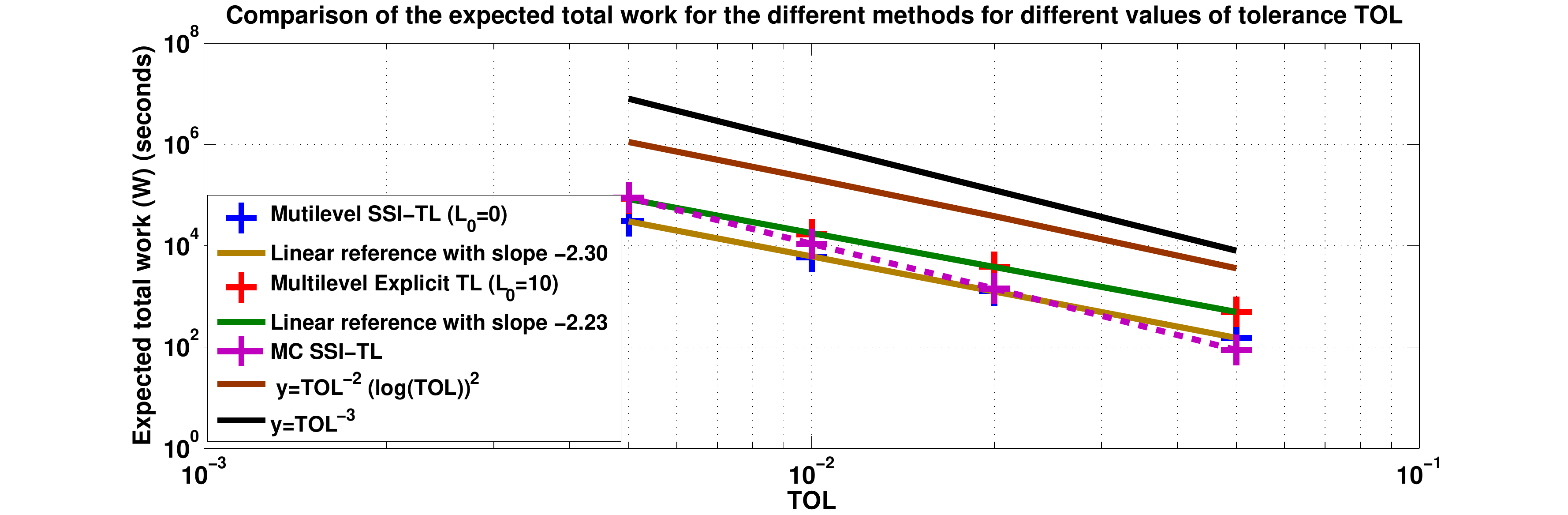}
\caption{Comparison of the expected total work for the different methods with different values of tolerance ($TOL$) for Example 1 using $100$ multilevel runs. The computational work, $W$, of both multilevel \name  and  multilevel \nameexp methods is of  $ \Ordo{TOL^{-2} (\log(TOL))^{2}}$ compared to $\Ordo{TOL^{-3}}$ for the MC \name method.}
\label{fig:Comparison of the expected total work for the different method}
 \end{center}
\end{figure}

To compare the computational work of the different methods, $100$ sets of multilevel calculations were performed for each value of tolerance, $TOL$. The actual work (runtime) was obtained using a 12-core Intel GLNXA64 architecture and MATLAB version R2014a. The results, given in  Table \ref{Comparison of the expected total work for the different methods for the decaying-dimerizing example}, indicate that we achieve the lowest computational cost using the multilevel \name estimator, which outperforms by about three times both the  multilevel \nameexp estimator and the MC \name estimator in terms of computational work. In addition, Figure \ref{fig:Comparison of the expected total work for the different method} shows that we achieve computational work, $W$, of $ \Ordo{TOL^{-2} (\log(TOL))^{2}}$ for both multilevel \name and  multilevel \nameexp  methods. This result confirms the computational advantage of MLMC over MC.
\begin{table}[H]
\centering
\begin{tabular}{l | r | r | r | r}
Method / TOL & $0.05$ &  $0.02$ &  $0.01$ & $0.005 $  \\
\hline \hline
 Multilevel \nameexp  &      500 (7)  &   3900 (21)      &  1.7+e04 (80)    &  8.6+e04 (570) \\
Multilevel \name   &  150 (3)  & 1300 (15)    &  5.9+e03 (47) &  3.1+e04 (265)  \\
 MC \name & 87 (1)  &  1400 (16) & 1.1+e04 (72) &  8.9+e04 (589)  \\
\hline
$W_{\text{MLMC}}^{\text{exp}}/W_{\text{MLMC}}^{\text{SSI}}$ & 3.30 &   2.97    &  2.85 &   2.80    \\
\hline
$W_{\text{MC}}^{\text{SSI}}/W_{\text{MLMC}}^{\text{SSI}}$ &      0.58    &  1.10   &    1.84   &  2.91 \\
\hline
\end{tabular}
\caption{Comparison of the expected total work for the different methods (in seconds) using $100$ multilevel runs for Example 1 $(\Ls=10  \; \text{and} \;  \Lc=0)$. The quantity in parentheses, \ie, $(\cdot)$,  refers to the standard deviation. The speedup factors given in the last two rows show that the multilevel \name  performs faster than the  multilevel \nameexp and the MC \name  methods.  }
\label{Comparison of the expected total work for the different methods for the decaying-dimerizing example}
\end{table}

\begin{figure}[H]
 \begin{center}
\includegraphics[scale=0.35]{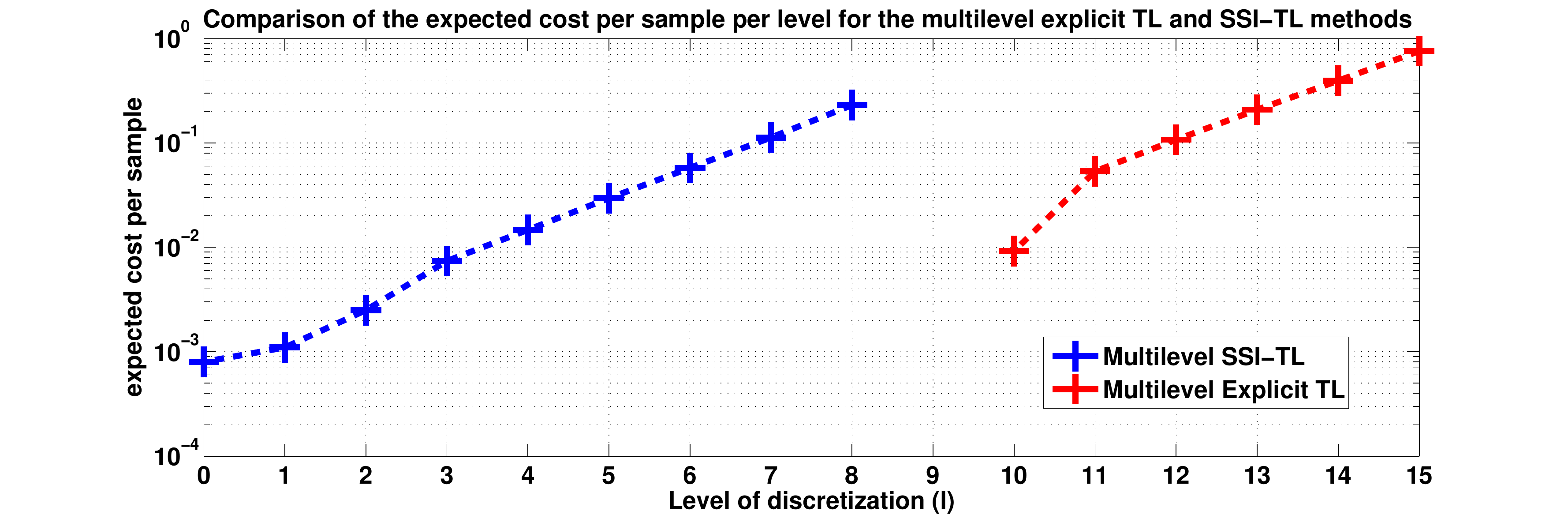}
\caption{Comparison of the expected cost per sample per level of the different methods for Example 1 using $10^5$ samples. The first observation corresponds to the time of a single path for the coarsest level and the other observations correspond to the time of the coupled paths per level. The computational gain due to using coarse-time discretization for \name  deteriorates due to the cost of Newton iterations as level $\ell$ increases. We note that we use a fixed number of Newton iterations equal to 3.}
\label{fig:Comparison of the expected cost per sample1}
 \end{center}
\end{figure}

This computational gain by the  multilevel \name  method  over the  multilevel \nameexp method is due to the lower cost of constructing single and coupled paths for the  \name method compared with the explicit one as illustrated by Figure \ref{fig:Comparison of the expected cost per sample1}. This figure shows that this computational gain, due to using coarse time discretization in the \name method, deteriorates due to the cost of Newton iterations. This is illustrated by having the same computational cost to generate coupled paths of \nameexp and \name while using different time step sizes (comparing the time to generate the coupled paths corresponding to level $\ell=8$ for the \name and the paths corresponding to level $\ell=13$ for the \nameexp). 
This observation motivated our idea of switching from the multilevel \name to the multilevel \nameexp method,  particularly when the value of $TOL$ is very small, implying a large value of the finest level of the drift-implicit MLMC, $L^{\text{imp}}$, such that $L^{\text{imp}} \geq \Ls $ ($\Ls$ refers to the first level of discretization from which the multilevel \nameexp method becomes stable).

\begin{rem}
 In future work, we will implement the method entirely in C++ to confirm the computational work gains.
\end{rem}

\begin{figure}[H]
\centering
\begin{minipage}{0.49\textwidth}
\includegraphics[scale=0.17]{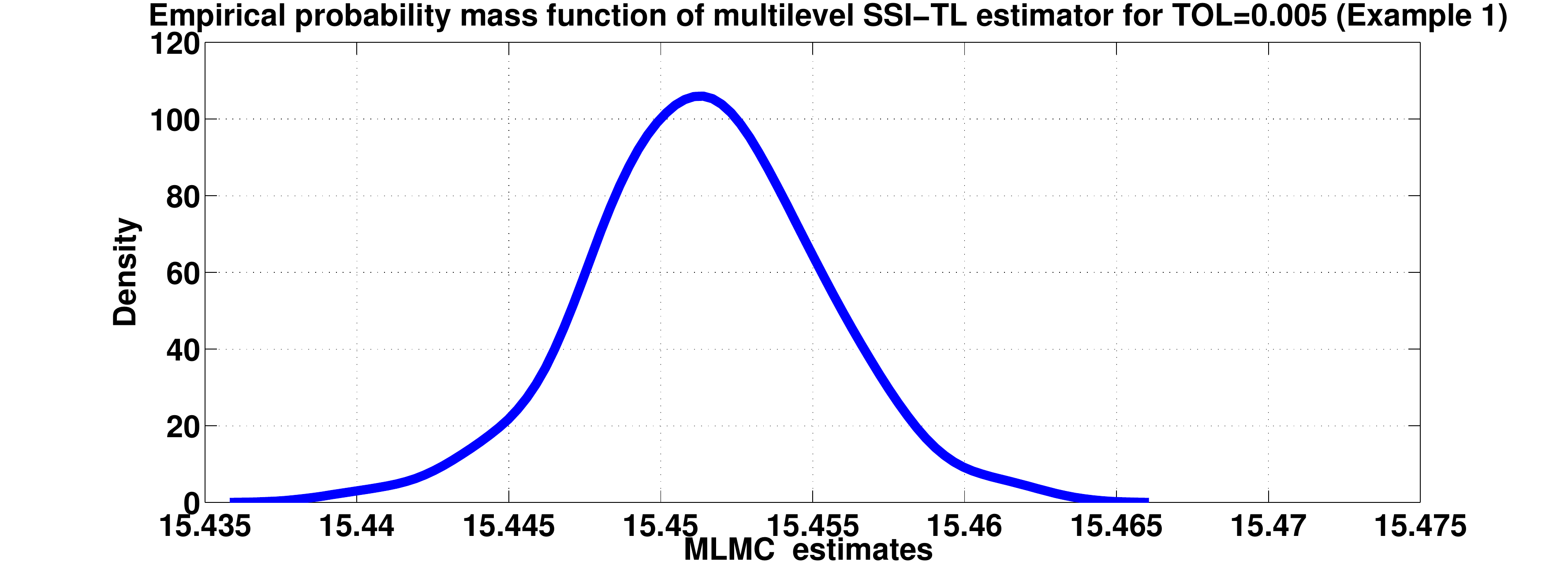}
\end{minipage}
\hfill
\begin{minipage}{0.49\textwidth}
\includegraphics[scale=0.18]{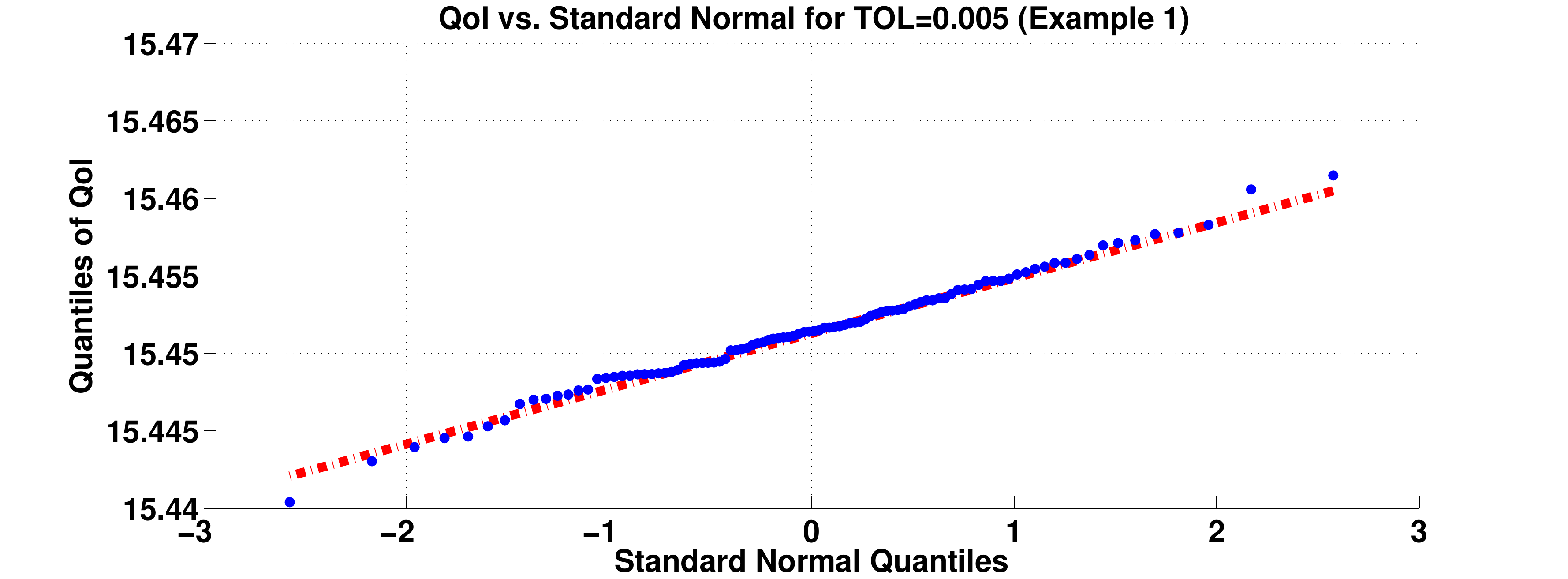}
\end{minipage}
\caption{Left: Empirical probability mass function for $100$ multilevel \name  estimates. Right: QQ-plot for the multilevel \name  estimates in Example 1.}
\label{fig:normality_test}
\end{figure}

\begin{figure}[H]
 \begin{center}
\includegraphics[scale=0.35]{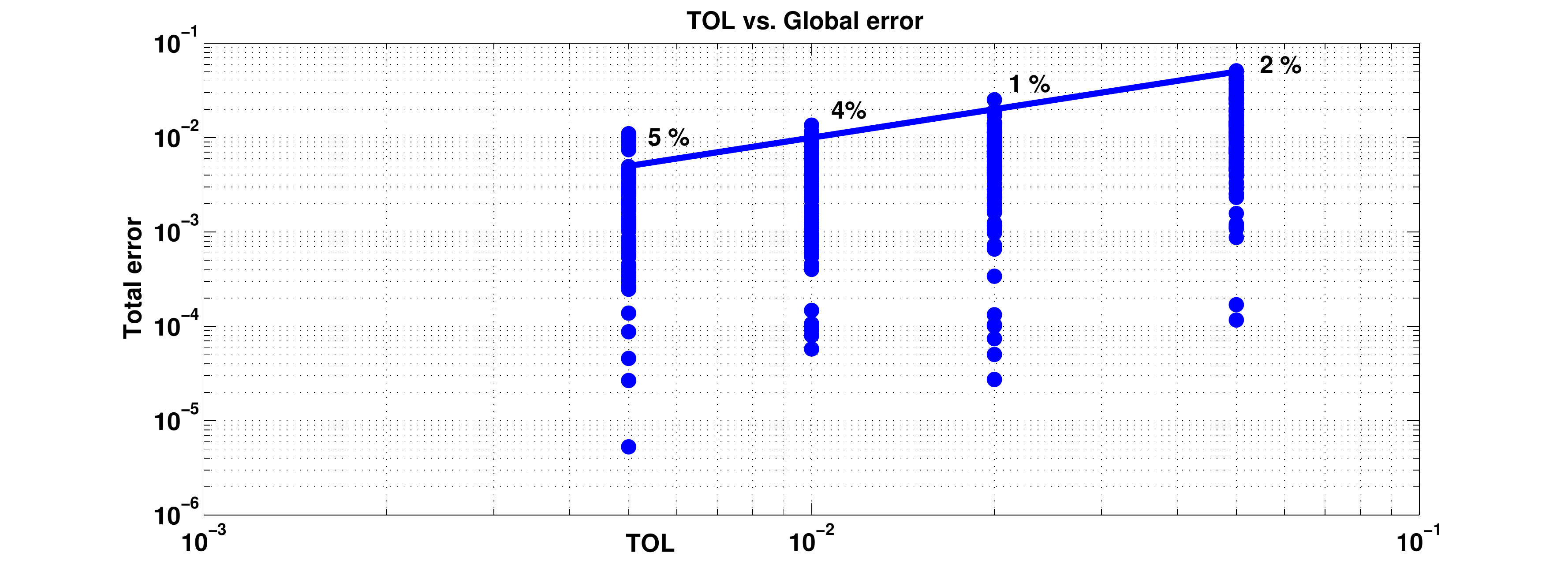}
\caption{TOL versus the actual computational error for Example 1. The numbers above the straight line show the percentage of runs that had errors larger than the required tolerance. We observe that in all cases, except for the smallest tolerance, the computational error follows the imposed tolerance with the expected confidence of $95 \%.$}
\label{fig:tol_global_exp1}
 \end{center}
\end{figure}

The QQ-plot and probability mass function plot in Figure \ref{fig:normality_test} show, for the smallest considered TOL, $100$ independent realizations of the multilevel \name  estimator. These plots, complemented by a Shapiro-Wilk normality test, validate our assumption about the Gaussian distribution of the statistical
error. In Figure (\ref{fig:tol_global_exp1}), we  show TOL versus the actual computational error. The prescribed tolerance is achieved with the required confidence of $95 \%$.

\subsubsection{Multilevel Hybrid \name Results ($L \geq \Ls$)}
\begin{figure}[H]
 \begin{center}
\includegraphics[scale=0.35]{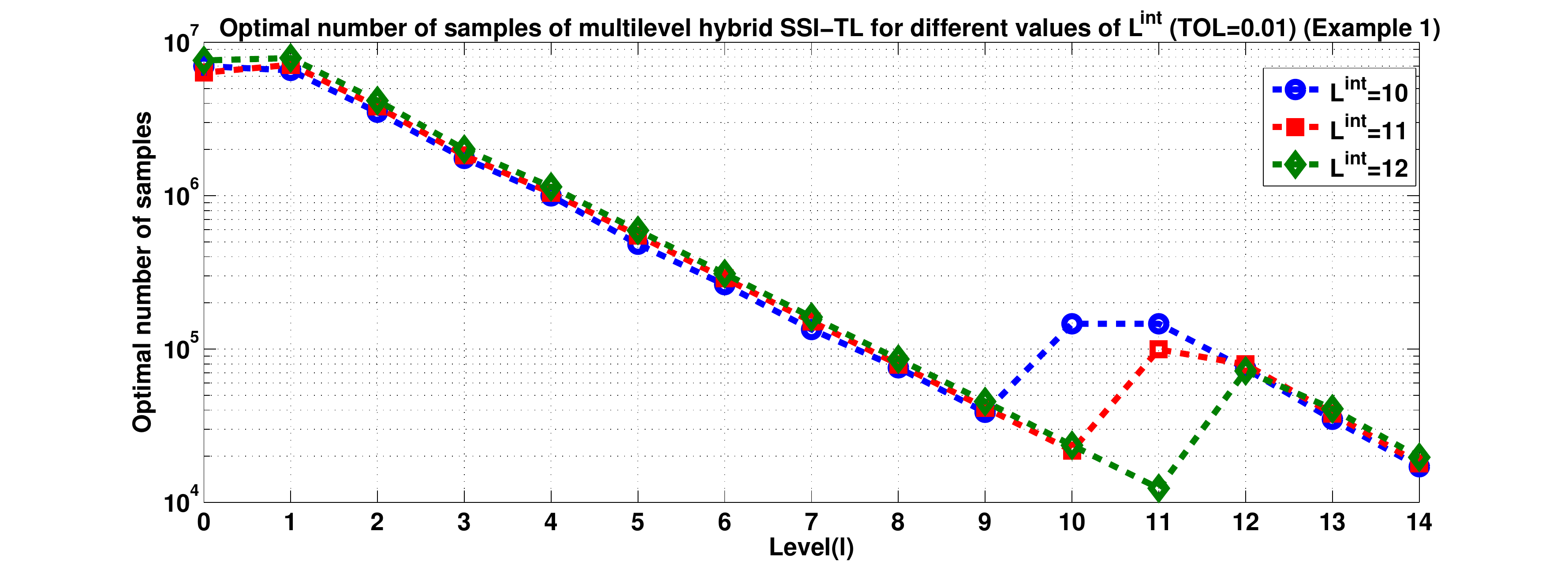}
\caption{Optimal number of samples of multilevel hybrid \name estimator for different values of $\Li$ ($TOL=0.01$) (Example 1). We observe jumps at the interface level, $\Li$. }
\label{fig:samples_hybrid_exp1_tol001}
 \end{center}
\end{figure}

\begin{figure}[H]
 \begin{center}
\includegraphics[scale=0.35]{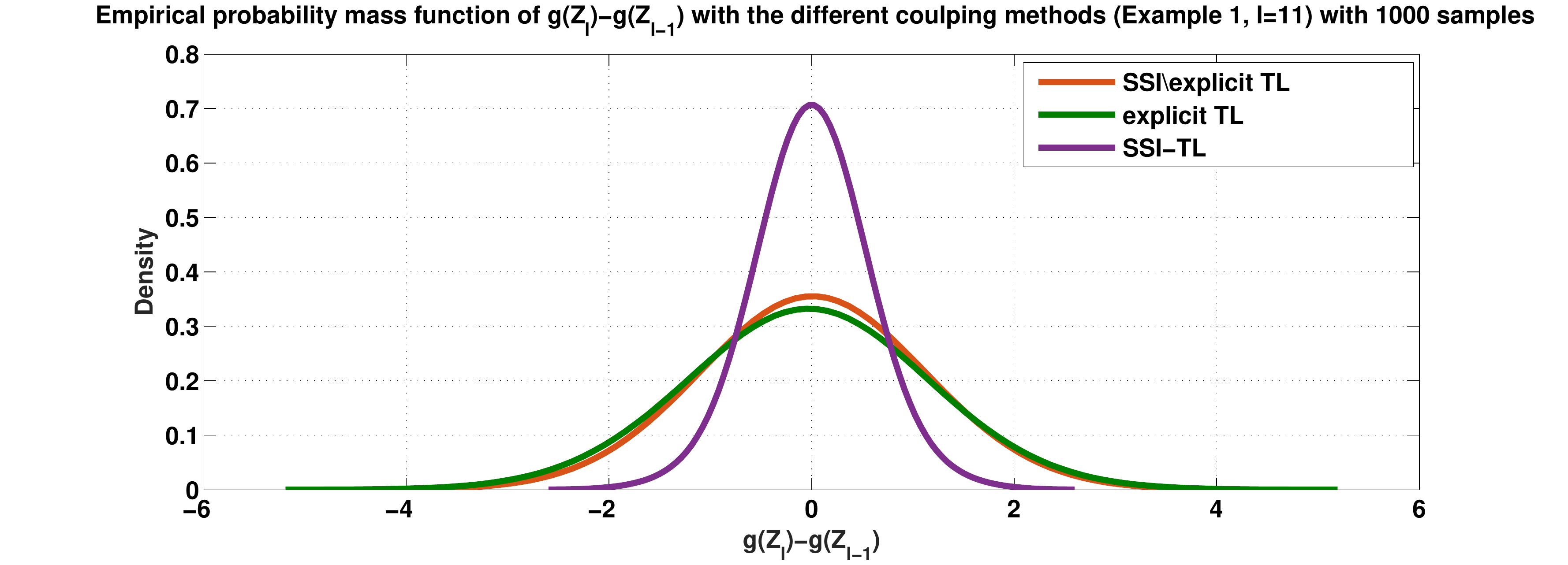}
\caption{Empirical probability mass function of $g(\mathbf{Z}_{\ell})-g(\mathbf{Z}_{\ell-1})$ simulated by the different coupling methods (Example 1, $\ell=11$) with $10^3$ samples.  The variance of the difference between coupled paths simulated by the \name method is smaller than the variance simulated either by the explicit or the coupling of the \name with the \nameexp method.}
\label{fig:density_comparison_hybrid_exp1_tol001}
 \end{center}
\end{figure}

\begin{figure}[H]
 \begin{center}
\includegraphics[scale=0.35]{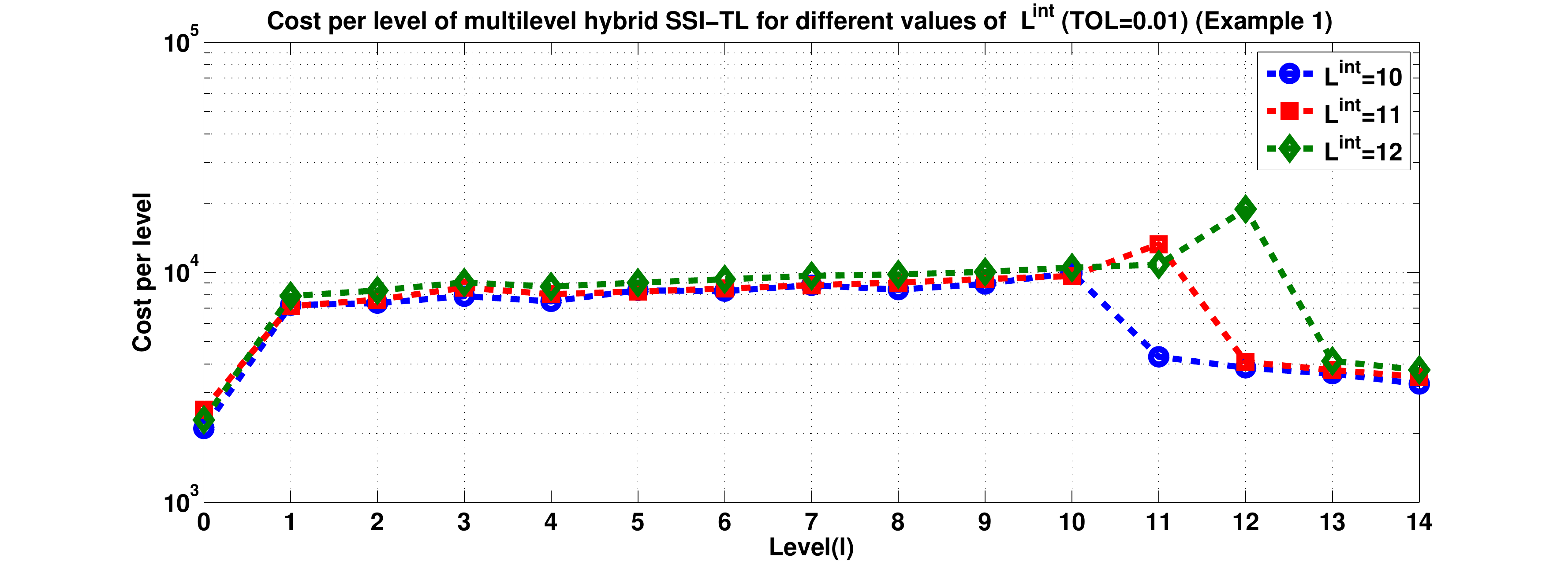}
\caption{Cost per level of the multilevel hybrid \name estimator for different values of $\Li$ ($TOL=0.01$) (Example 1). This figure explains why the lowest computational cost is most likely to be achieved for $\Li = \Ls$ (see analysis in Section \ref{sec:Hybrid Unbiased Tau-leap MLMC Estimator}).}
\label{fig:cost_hybrid_exp1_tol001}
 \end{center}
\end{figure}

Here, we show results for our multilevel hybrid \name estimator for different values of the interface level $\Li$. 
Table \ref{Comparison of the expected total work for the different hybrid methods for example 1} shows that we achieve the lowest computational cost for $\Li=\Ls$. This can be explained by  analyzing the cost per level of the multilevel hybrid \name estimator (see Figure \ref{fig:cost_hybrid_exp1_tol001} and using relation (\ref{hybrid_work_analysis})). Table \ref{Comparison of the expected total work for the different hybrid methods for example 1}  also shows that our estimator outperforms the explicit one by about three times. This gain can be more important for very small values of tolerance, $TOL$. Figure \ref{fig:samples_hybrid_exp1_tol001} shows the optimal number of samples in our \name setting. We observe jumps at the interface level, $\Li$, which are a consequence of jumps in the variance of the coupled paths. 
This is due to the fact that the variance of coupled paths, for the same discretization level, simulated by the  \name is smaller than the variance simulated either by the explicit or the coupling of the \name with the \nameexp method (see Figure \ref{fig:density_comparison_hybrid_exp1_tol001}). 

\begin{table}[H]
\centering
\begin{tabular}{l | r | r }
Method / TOL & $0.01$ &  $0.005$   \\
\hline \hline
 Multilevel \nameexp $(\Ls=10)$   &     1.7e+04(80)    &  8.6e+04 (570)  \\
Multilevel Hybrid \name  $(\Lc=0,\Li=10) $  & 7.2e+03 (52)  & 2.8e+04 (239) \\
Multilevel Hybrid \name  $(\Lc=0,\Li=11) $  & 8.7e+03 (55) &    3.2e+04 (271)   \\
Multilevel Hybrid \name $(\Lc=0,\Li=12) $  &  1.1e+04 (69) & 3.9e+04 (246)     \\
\hline
$W_{\text{MLMC}}^{\text{exp}}/W_{\text{MLMC}}^{\text{hyb}}(\Li=10)$ &  2.36 &   3.07      \\
\hline
\end{tabular}
\vspace{0.1cm}
\caption{Comparison of the expected total work for the different methods (in seconds) using $100$ multilevel runs for Example 1. The quantity in parentheses refers to the standard deviation. The optimal value of the interface level is $\Li=\Ls$. 
The speedup factor, represented by the last row, can reach a factor of $3$.  }
\label{Comparison of the expected total work for the different hybrid methods for example 1}
\end{table}

\subsection{Example 2}
This example was studied in \cite{Gillespie03} and is given by the following reaction set
\begin{align}\label{example2}
S_{1} &\overset {c_1}{\rightarrow} S_{3} ,\nonumber \\
S_{3} &\overset {c_2}{\rightarrow} S_1  \nonumber \\
S_{1}+S_{2} &\overset {c_3}{\rightarrow} S_{1}+S_{4}.
\end{align}
Since the total number of $S_1$ and $S_3$ molecules is constant (say $K$), and since we can ignore the by-product, $S_4$, this system can be represented by three reactions and two variables, $\mathbf{X}=(X_1 ,X_2)$, which are numbers of $S_1$ and $S_2$ molecules, respectively. The stoichiometric vectors are $\boldsymbol{\nu}_{1}=(-1,0)^{T}$, $\boldsymbol{\nu}_{2}=(1,0)^{T}$, and $\boldsymbol{\nu}_{3}=(0,-1)^{T}$ and  the corresponding propensity functions are
\begin{equation}
a_{1}(\mathbf{X}) = c_{1} X_{1},\: a_{2}(\mathbf{X}) = c_{2} (K - X_{1}),\: a_{3}(\mathbf{X})  = c_{3} X_{1} X_2.
\end{equation}
We chose the same values for the different parameters as in \cite{Gillespie03}: $c_1=c_2=10^{5}$, $c_3=5 \times 10^{-3}$, $K=2 \times 10^{4}$ and  initial condition $\mathbf{X}(0)=(10^{4},10^{2})^{T}$. This setting implies that the stability limit of the \nameexp is $\tau_{\text{exp}}^{\text{lim}} \approx  10^{-5}$.  We consider the final time, $T= 0.01$ seconds. In the following numerical experiments, we are interested in approximating $\expt{X_{2}(T)}$.

\subsubsection{Multilevel  \name Results ($L \leq \Ls$)}
Similarly to what we did in the first example, we checked whether the coupling procedure is reliable in terms of convergence properties and we followed the same steps shown in Example 1. Now, we want to compare the performance of our multilevel hybrid \name algorithm, reduced to the multilevel \name algorithm   when $L \leq \Ls$, to the performance of the  multilevel \nameexp and  MC \name methods. Here, we fix $\Ls=11$ for the  multilevel \nameexp  estimator to ensure the numerical stability ($2^{-\Ls} T < 10^{-5}$). We also vary the prescribed tolerance, $TOL$, to investigate its effect on the performance of the tested methods. We note that, in this example, we use a fixed number of Newton iterations equal to 3 for the \name method.

\begin{figure}[H]
 \begin{center}
\includegraphics[scale=0.35]{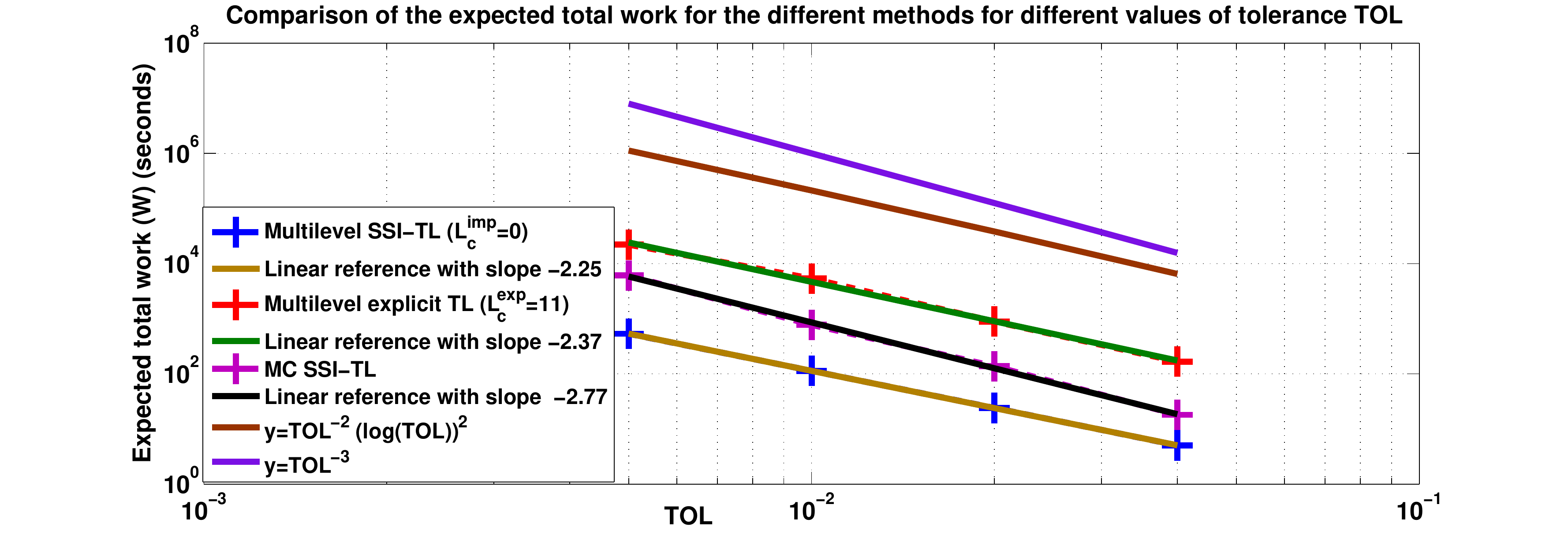}
\caption{Comparison of the expected total work for the different methods with different values of tolerance ($TOL$) for Example 2 using $100$ multilevel runs. The computational work, $W$, of both multilevel \name and multilevel \nameexp methods is   $ \Ordo{TOL^{-2} (\log(TOL))^{2}}$ compared to $\Ordo{TOL^{-3}}$ for the MC \name.}
\label{fig:Comparison of the expected total work for the different method_exp2}
 \end{center}
\end{figure}

\begin{figure}[H]
 \begin{center}
\includegraphics[scale=0.35]{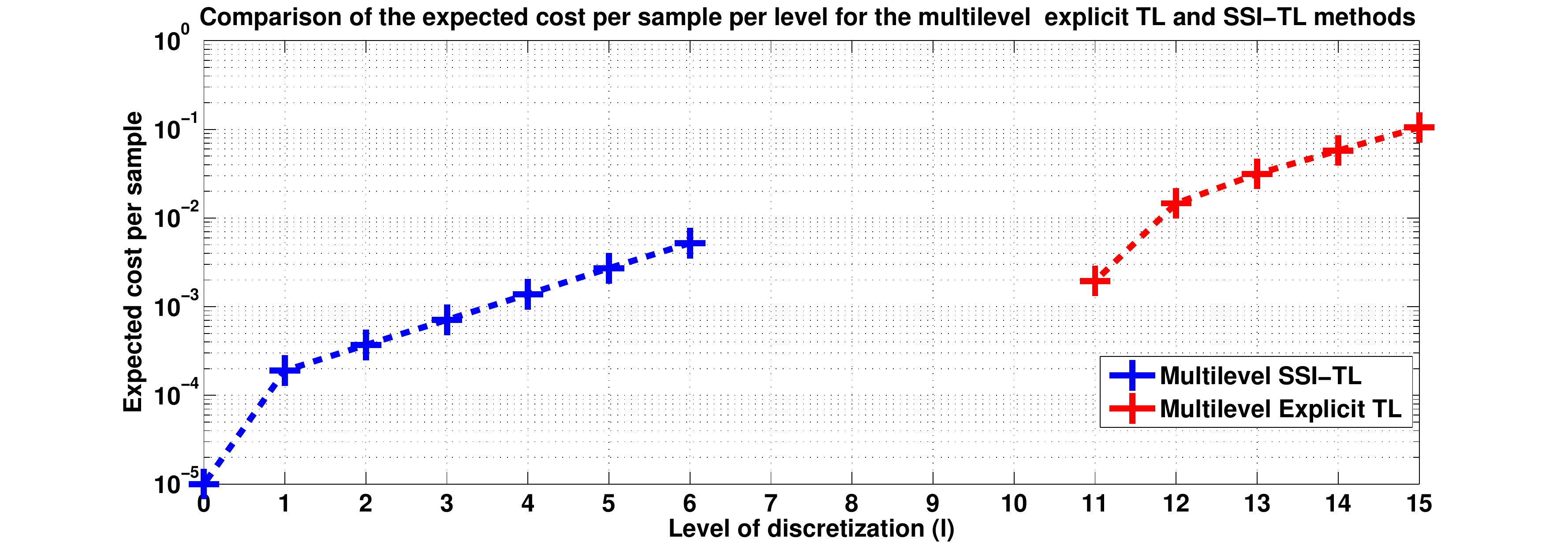}
\caption [Comparison of the expected cost per sample per level for the different methods for Example 2 using $10000$ samples.] {Comparison of the expected cost per sample per level for the different methods for Example 2 using $10000$ samples. The first observation corresponds to the time of a single path for the coarsest level and the other observations correspond to the time of the coupled paths per level.}
\label{fig:Comparison of the expected cost per sample2}
 \end{center}
\end{figure}

\begin{figure}[H]
\centering
\begin{minipage}{0.49\textwidth}
\includegraphics[scale=0.17]{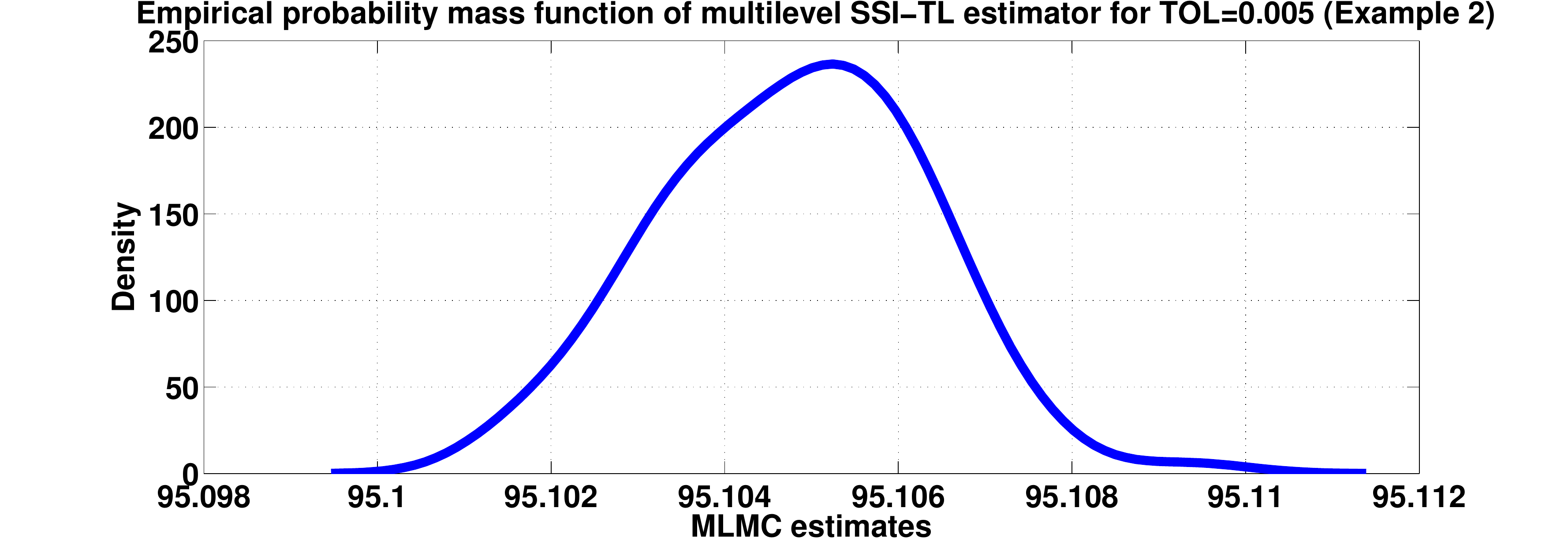}
\end{minipage}
\hfill
\begin{minipage}{0.49\textwidth}
\includegraphics[scale=0.18]{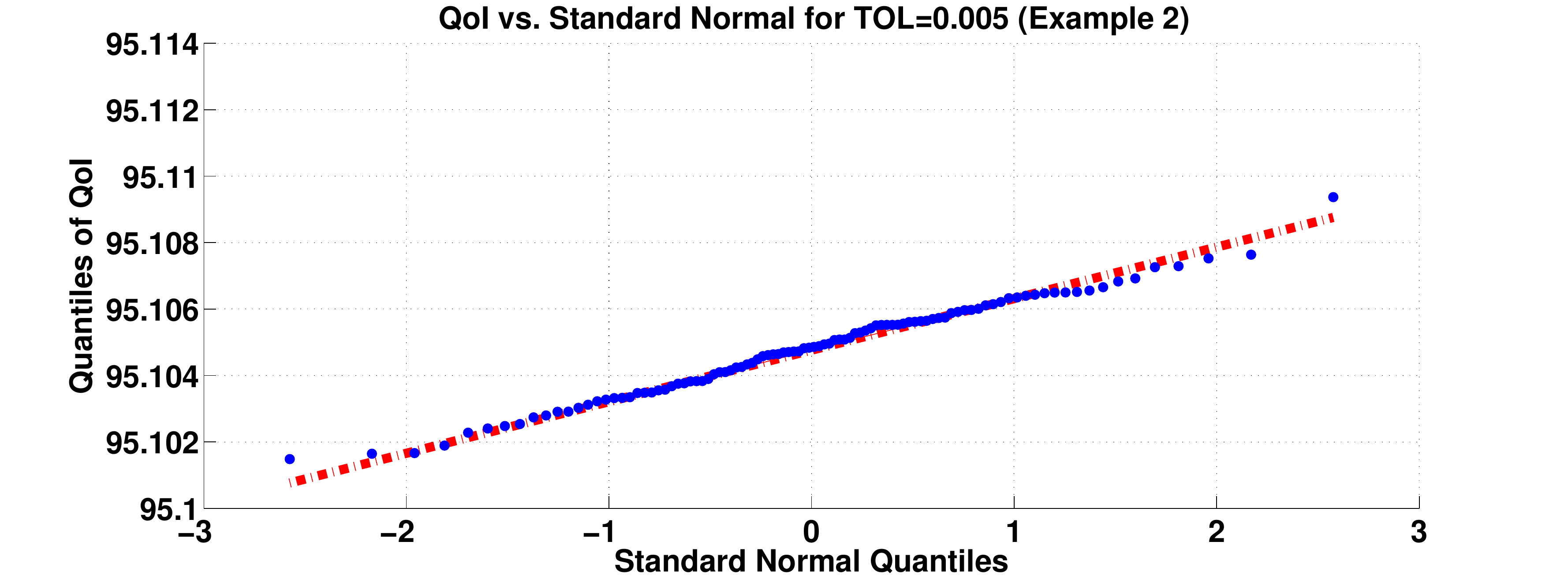}
\end{minipage}
\label{fig:normality_test2}
\caption{Left: Empirical probability mass function for $100$ multilevel  \name estimates. Right: QQ-plot for the multilevel \name estimates in Example 2. }
\end{figure}

\begin{figure}[H]
 \begin{center}
\includegraphics[scale=0.35]{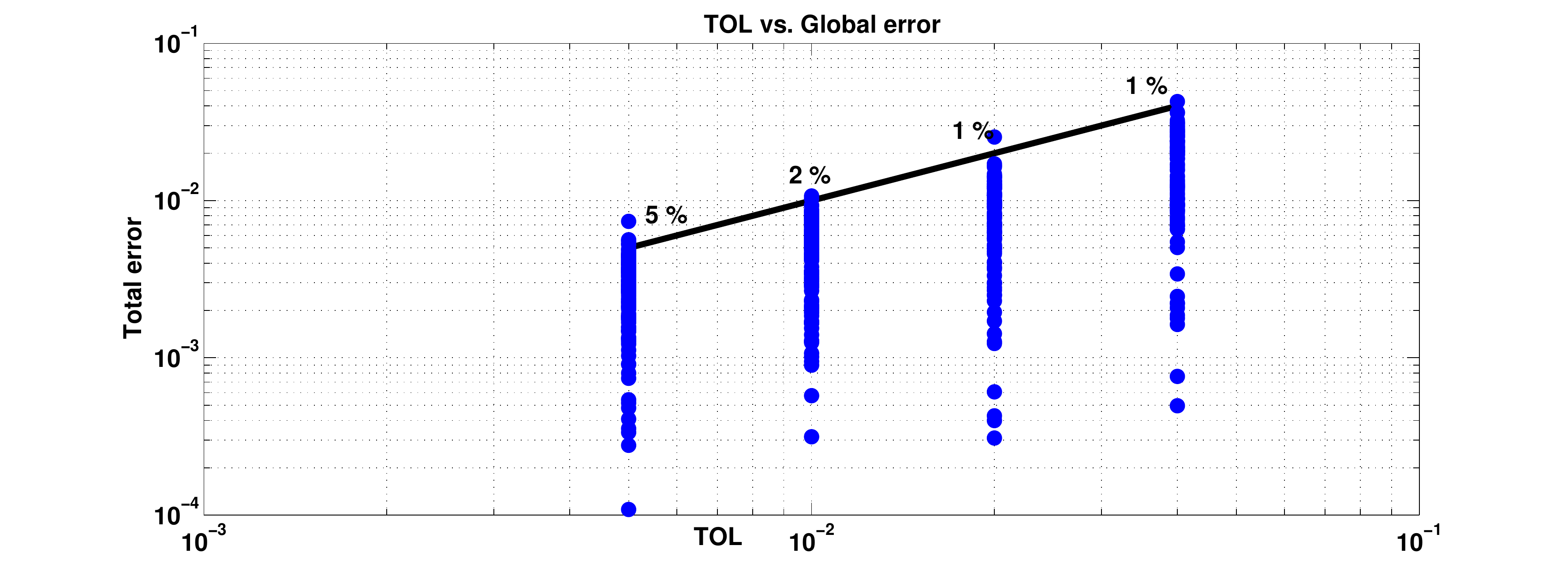}
\caption{TOL versus the actual computational error for Example 2. The numbers above the straight line show the percentage of runs that had errors larger than the required tolerance. We observe that in all cases, except for the smallest tolerance, the computational error follows the imposed tolerance with the expected confidence of $95 \%.$}
\label{fig:tol_global_exp2}
 \end{center}
\end{figure}

To compare the computational work of the different methods, $100$ sets of multilevel calculations were performed for each value of tolerance $TOL$. The results presented in Table \ref{Comparison of the expected total work for the different methods for example 2} and Figure \ref{fig:Comparison of the expected total work for the different method_exp2} indicate that we achieve the lowest computational cost with the multilevel \name estimator, which outperforms by about $40$ times the  multilevel \nameexp estimator and by  about $11$ times the  MC \name estimator in terms of computational work. Similarly to the first example, Figure \ref{fig:Comparison of the expected total work for the different method} shows that we achieve computational work, $W$, of $ \Ordo{TOL^{-2} (\log(TOL))^{2}}$ for both multilevel \name and  multilevel \nameexp methods. This result confirms the computational advantage of MLMC over MC.
\begin{table}[H]
\centering
\begin{tabular}{l | r | r | r | r}
Method / TOL & $0.04$ &  $0.02$ &  $0.01$ & $0.005 $  \\
\hline \hline
 Multilevel \nameexp &   170   (4)  &   890   (9)  &     5300 (45) &   2.2e+04 (96) \\
Multilevel \name   &   5 (0.2)   &   240  (0.8)     & 110 (3)  &    5.3e+02  (7) \\
 MC \name &  18 (0.6) & 140   (3) &  770 (8) &    6e+03 (48)  \\
\hline
$W_{\text{MLMC}}^{\text{exp}}/W_{\text{MLMC}}^{\text{SSI}}$ &   33.40  &    37.04   &  47.33  &       41.27\\
\hline
$W_{\text{MC}}^{\text{SSI}}/W_{\text{MLMC}}^{\text{SSI}}$  &    3.60 &   5.70    &   6.84   &  11.37   \\
\hline
\end{tabular}
\vspace{0.1cm}
\caption{Comparison of the expected total work for the different methods (in seconds) using $100$ multilevel runs for Example 2 $(\Ls=11 \: \: \text{and} \: \: \Lc=0)$. The quantity in parentheses  refers to the standard deviation. The speedup factors are more important in this example. From the last two rows, we notice that the multilevel  \name estimator outperforms by about $40$ times the multilevel \nameexp estimator and  by about $11$ times the MC \name estimator in terms of computational work.}
\label{Comparison of the expected total work for the different methods for example 2}
\end{table}
This computational gain of the multilevel \name method over the  multilevel \nameexp method is, as in the previous example, due to the lower cost of constructing single and coupled paths by \name compared to the cost by the explicit method as illustrated by Figure \ref{fig:Comparison of the expected cost per sample2}. Similarly to the first example, we can  notice from this figure that the computational gain due to using coarse-time discretization for the \name method deteriorates due to the cost of Newton iterations.
The QQ-plot and probability mass function plot in Figure \ref{fig:normality_test2} show, for the smallest considered TOL, $100$ independent realizations of the multilevel \name estimator. These plots, complemented with a Shapiro-Wilk normality test, validates our assumption about the Gaussian distribution of the statistical
error. In Figure \ref{fig:tol_global_exp2}, we  show TOL versus the actual computational error. The prescribed tolerance is achieved with the required confidence of $95 \%$.

\subsubsection{Multilevel Hybrid \name Results ($L \geq \Ls$)}
\begin{figure}[H]
 \begin{center}
\includegraphics[scale=0.35]{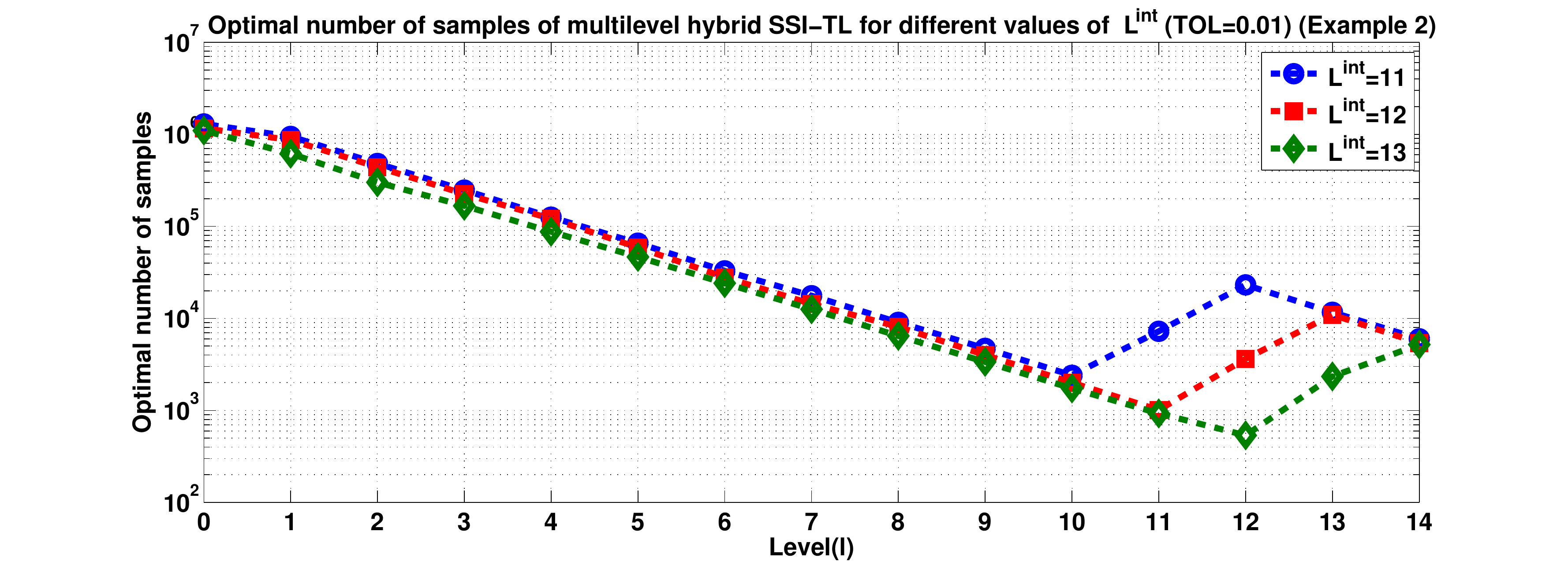}
\caption{Optimal number of samples of multilevel hybrid \name for different values of $\Li$ ($TOL=0.01$) (Example 2). We observe jumps at the interface level, $\Li$.}
\label{fig:samples_hybrid_exp2_tol001}
 \end{center}
\end{figure}

\begin{figure}[H]
 \begin{center}
\includegraphics[scale=0.35]{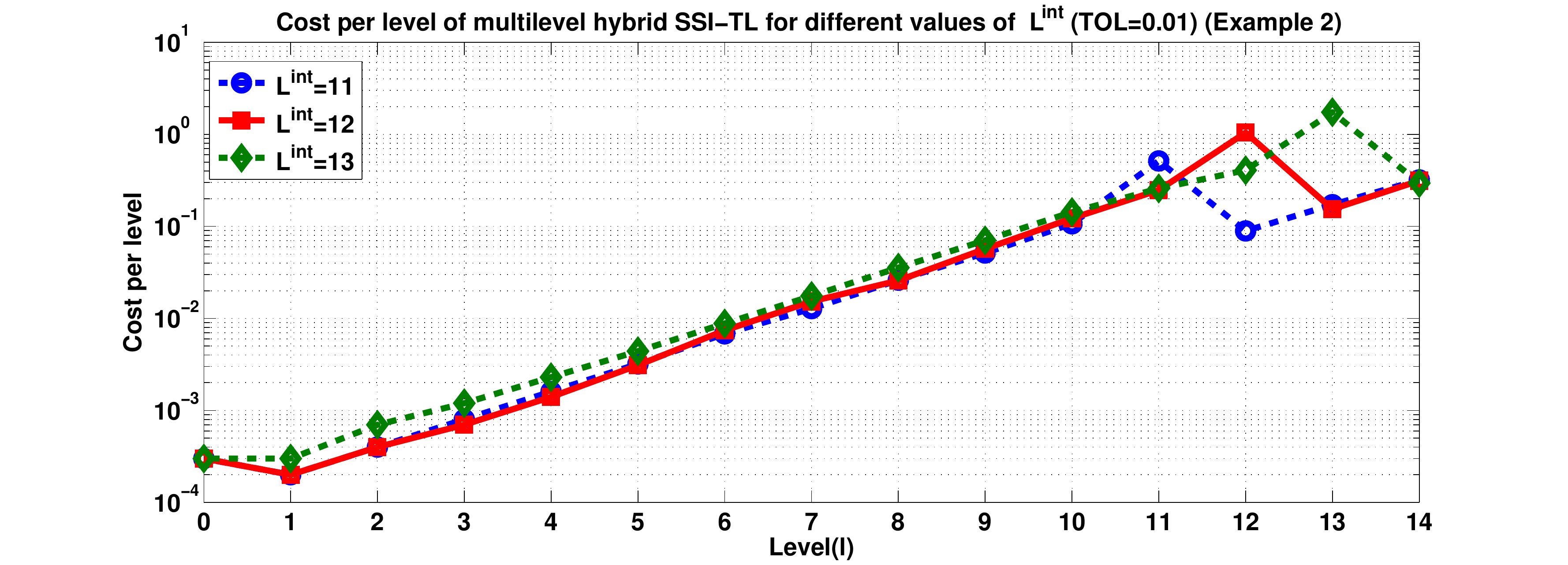}
\caption{Cost per level of the multilevel hybrid \name  estimator for different values of $\Li$ ($TOL=0.01$) (Example 2). This figure explains why the lowest computational cost is most likely to be achieved for $\Li = \Ls$ (see analysis  in Section \ref{sec:Hybrid Unbiased Tau-leap MLMC Estimator}).}
\label{fig:cost_hybrid_exp2_tol001}
 \end{center}
\end{figure}

Here, we present the results of the multilevel hybrid \name estimator for different values of the interface level $\Li$. Table \ref{Comparison of the expected total work for the different hybrid methods for example 2} shows that we achieve the lowest computational cost for, $\Li=\Ls$. This can be explained by  analyzing the cost per level of the multilevel hybrid \name estimator (see Figure \ref{fig:cost_hybrid_exp2_tol001} and using relation (\ref{hybrid_work_analysis})). Table \ref{Comparison of the expected total work for the different hybrid methods for example 2} also shows that our multilevel estimator  outperforms the explicit one by  about seven times. This gain can be more important for very small values of tolerance, $TOL$. Figure \ref{fig:samples_hybrid_exp2_tol001} shows the optimal number of samples in the \name setting. We observe jumps at the interface level, $\Li$, which are a consequence of jumps in the variance of the coupled paths.  
\begin{table}[H]
\centering
\begin{tabular}{l | r | r }
Method / TOL & $0.01$ &  $0.005$   \\
\hline \hline
 Multilevel \nameexp $(\Ls=11)$   &     5.3e+03 (45) &   2.2e+04 (96)  \\
Multilevel Hybrid \name $(\Lc=0,\Li=11) $  &  9.4e+02 (11)  & 3.2e+03 (15)  \\
Multilevel Hybrid \name  $(\Lc=0,\Li=12) $  & 9.9e+02 (12)  &   3.3e+03 (15)    \\
Multilevel Hybrid \name   $(\Lc=0,\Li=13) $  & 1.3e+03 (13)  &   3.6e+03 (19)   \\
\hline
$W_{\text{MLMC}}^{\text{exp}}/W_{\text{MLMC}}^{\text{hyb}} (\Li=11)$ &  5.64 &     6.87      \\
\hline
\end{tabular}
\vspace{0.1cm}
\caption{Comparison of the expected total work for the different methods (in seconds) using $100$ multilevel runs for Example 2. The quantity in parentheses refers to the standard deviation.  The optimal value of the interface level is $\Li=\Ls$. The speedup factor, given in the last row, can reach a factor of $7$.}
\label{Comparison of the expected total work for the different hybrid methods for example 2}
\end{table}

\begin{rem}
As we argued in the introduction, from Figure \ref{fig:Comparison with tau Rock}, we check that simulating single paths with the \name method is more computationally efficient than with the $\tau$-ROCK method, especially in the case of large time steps.

\begin{figure}[H]
 \begin{center}
\includegraphics[scale=0.35]{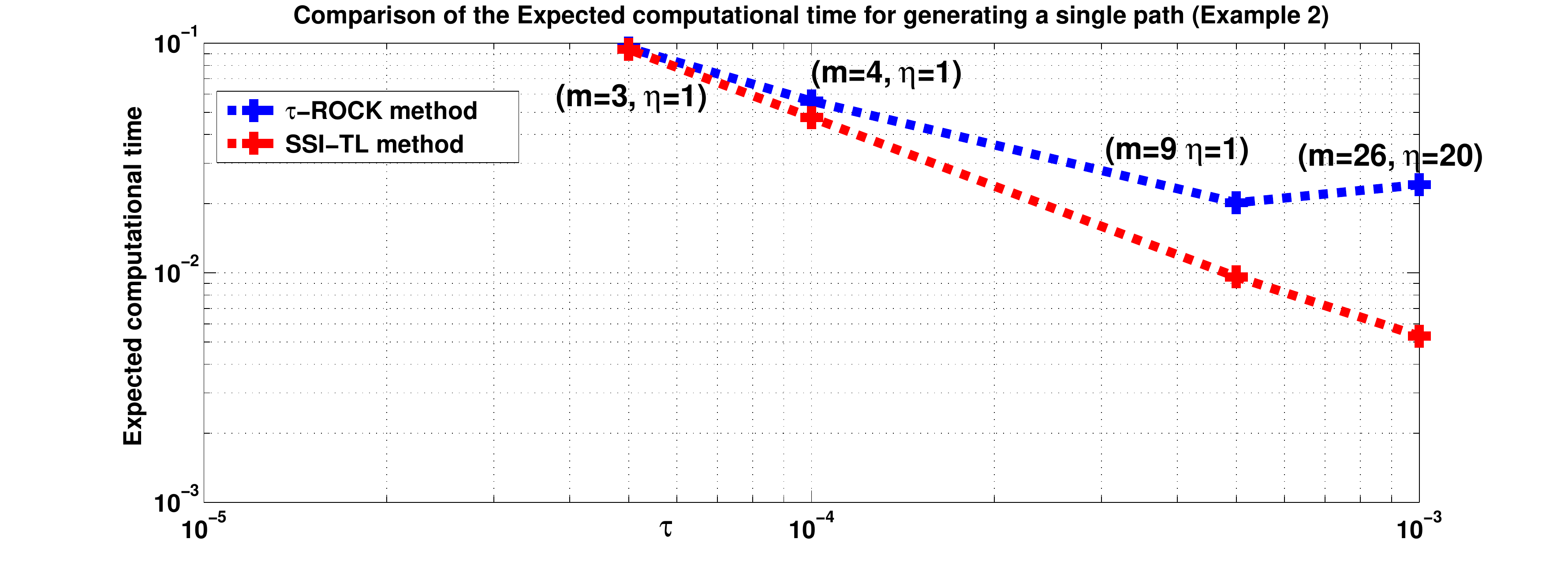}
\caption{Comparison between the \name and the $\tau$-ROCK methods in terms of the expected computational time, using $10^5$ samples, for generating a single path (Example 2 in Section \ref{sec:examples}). The parameters $m$ ( the stage number that controls the stability region) and $\eta$ (the damping parameter that controls the damping effect of the variance) are chosen in such a way that we obtain nearly optimal work for the $\tau$-ROCK method except for $\tau=10^{-3}$. For the sake of comparison, we also show the same values given in \cite{abdulle2010chebyshev}}.
\label{fig:Comparison with tau Rock}
  \end{center}
\end{figure}
\end{rem}

 
\section{Conclusions and Future Work}
\label{conclusions}
In this work,  we adapt the idea of split-step backward-Euler \cite{higham2002strong} for Brownian SDEs to SDEs driven by Poisson random measures. We extend the \name idea to the multilevel setting and introduce the multilevel hybrid \name estimator with the aim of reducing the computational work needed to produce an estimate of $\expt{g(\mathbf{X}(T))}$, within a fixed tolerance, $TOL$,  with a given level of confidence. Our estimator couples \name paths at the coarser discretization levels until a certain interface level, $\Li$, is reached. The $\Li$ level can be reached or not, depending on user-selected tolerance requirements. When the required number of levels is greater than $\Li$, our estimator couples \nameexp paths at the finer levels and can be even become unbiased by coupling the finest level with a modified next reaction method \cite{anderson:214107}.
 
Our proposed estimator is useful in systems with the presence of slow and fast timescales (stiff systems). In such situations, the multilevel estimator given in \cite{Anderson2012} is not computationally efficient due to  numerical stability constraints. Through our numerical experiments, we obtained substantial gains with respect to both the  multilevel \nameexp and the single-level \name methods. We also showed, that for large values of $TOL$, the multilevel hybrid \name method has the same order of computational work as does the  multilevel \nameexp method, which is of  $\Ordo{TOL^{-2} \log(TOL)^{2}}$  \cite{Anderson_Complexity}, but with a smaller constant.
{Although it is possible to obtain a computational work of  $\Ordo{TOL^{-2}}$ by coupling with a pathwise exact path at the deepest level,  this is of little help in the class of stiff problems due to the high number of exact steps}.	
{
 As we stressed in Section \ref{sec:Hybrid Unbiased Tau-leap MLMC Estimator}, one potential direction of future research is applying the CMLMC algorithm, introduced in \cite{collier_continuation_MLMC}, in our context. In fact, it may provide an optimal split between bias and statistical errors, implying an improvement in the the computational complexity of our MLMC estimator. Future extensions may also involve  hybridization techniques involving methods that deal with
non-negativity of species (see \cite{alvaro_ML}). These hybridization techniques can be improved by using the idea of adaptivity introduced in \cite{tau_control_variate,lester2015adaptive}, which allows us to construct adaptive hybrid multilevel estimators by switching between \name, \nameexp and exact SSA within the course of a single sample. We also intend to extend the $\tau$-ROCK method \cite{abdulle2010chebyshev} to the multilevel setting and compare it with the multilevel hybrid \name method. The main challenge of this task is to couple two consecutive paths based on the $\tau$-ROCK method. Finally, our techniques can be extended to the context of spatial inhomogeneities described, for instance, by graphs and/or continuum volumes.}

\section*{Acknowledgments}
Research reported in this publication was supported by competitive research funding from King Abdullah University of Science and Technology (KAUST). C. Ben Hammouda, A. Moraes and R. Tempone are members of the KAUST SRI Center for Uncertainty Quantification in the Computer, Electrical and Mathematical Science and Engineering Division, KAUST. We  would also like to thank the anonymous referees for the constructive
feedback that helped improve this manuscript.

\appendix
\section{Results of MLMC for Example 1}\label{RME1}

\begin{table}[H]
\centering
\begin{tabular}{l | r | r | r | r }
 Level / tol  & $0.05$ &  $0.02$ &  $0.01$ & $0.005 $   \\
\hline \hline
$15$ &  -  &  -    & - & $ 2.5  \times 10^{4} $\\
$14 $&  -  &  -     &  $3.6  \times 10^{4}$ &    $ 5.3 \times 10^{4}$ \\
$13$  & -  &   $ 1.8 \times 10^{4} $& $7.3  \times 10^{4}$ &  $3.1 \times 10^{5} $ \\
$12$ & $4.5 \times 10^{3}$ & $ 3.4 \times 10^{4}$  &$ 1.4  \times 10^{5}$  &  $ 5.3  \times 10^{5} $\\
$11$ & $9.6 \times 10^{3} $  & $ 7.9 \times 10^{4}$ &  $ 3.3 \times 10^{5}$  &  $ 1.3 \times 10^{6}$ \\
$10$ &  $1.6 \times 10^{5} $  & $ 1.2 \times 10^{6}$  & $5 .0 \times 10^{6} $   & $2.1  \times 10^{7}$  \\
\end{tabular}
\vspace{0.1cm}
\caption{Optimal number of samples per level for the  multilevel \nameexp for example 1 $(L_{0}=10)$}
\label{Samples per level for the MLMC explicit tau-leap for the decaying-dimerizing example}
\end{table}
    
\begin{table}[H]
\centering
\begin{tabular}{l | r | r | r | r }
 Level / tol  & $0.05$ &  $0.02$ &  $0.01$ & $0.005 $   \\
\hline \hline
$8$ &  - &   -     & - & $2.0 \times 10^{5}$\\
$7$ &  - &  -  & $8.7 \times 10^{4}$ & $3.9 \times 10^{5}$\\
$6$ & - &   $ 3.5 \times 10^{4} $   & $1.7 \times 10^{5}$& $ 7.3 \times 10^{5} $ \\
$5$ & $ 8.8 \times 10^{3} $ & $ 7.0 \times 10^{4} $ &$ 3.1 \times 10^{5} $&  $1.4 \times 10^{6} $\\
$4$ &   $1.7 \times 10^{4} $ &  $1.3 \times 10^{5}$ & $ 6.1 \times 10^{5} $ & $ 2.7 \times 10^{6}$ \\
$3$ &   $3.3 \times 10^{4} $ & $2.6 \times 10^{5}$  & $10^{6}$   & $5.2 \times 10^{6}$   \\
$2$ &  $5.8 \times 10^{4} $ &  $4.6 \times 10^{5} $ &  $1.9 \times 10^{6}$   &   $9.3 \times 10^{6} $ \\
$1$ & $ 1.2 \times 10^{5}$  &   $9.4 \times 10^{5}$& $3.9 \times 10^{6}$  &   $  1.9 \times 10^{7} $  \\
$0 $&  $6.7 \times 10^{5}$ & $3.5 \times 10^{6} $ &  $ 2.4 \times 10^{7}  $ &   $10^{8}$  \\
\end{tabular}
\vspace{0.1cm}
\caption{Optimal number of samples per level for the multilevel  \name for example 1 $(L_{0}=0)$}
\label{Samples per level for the drift implicit MLMC tau-leap for the decaying-dimerizing example}
\end{table}
    
\begin{table}[H]
\centering
\begin{tabular}{l | r | r | r | r }
 Level / tol  & $0.05$ &  $0.02$ &  $0.01$ & $0.005 $   \\
\hline \hline
$8$&  - &   -     & - & $1.9 \times 10^{5}$\\
$7$ &  - &  -  & $7.7 \times 10^{4}$ &  $3.6 \times 10^{5}$\\
$6$ & - &   $ 3.1 \times 10^{4} $   &  $1.5 \times 10^{5} $ & $6.7 \times 10^{5} $ \\
$5$ & $8.0 \times 10^{3}$  &  $6.1 \times 10^{4} $ & $ 2.9 \times 10^{5} $& $ 1.3 \times 10^{6}$ \\
$4$ &   $1.6 \times 10^{4}$  & $ 1.2 \times 10^{5} $& $5.5 \times 10^{5}$  &  $2.6 \times 10^{6}$ \\
$3$ &   $2.8 \times 10^{4}$  & $2.1 \times 10^{5} $ &   $10^{6}  $ & $ 4.6 \times 10^{6} $  \\
$2$ & $ 5.7 \times 10^{4} $ & $ 4.4 \times 10^{5}$  &  $2.0 \times 10^{6} $  &  $10^{7}$  \\
$1$ &  $1.1 \times 10^{5} $ &  $ 8.5 \times 10^{5} $& $3.9 \times 10^{6} $ &   $  1.9 \times 10^{7} $  \\
$0$ &  $1.1 \times 10^{5} $ & $ 8.1 \times 10^{5}$  &   $4.0 \times 10^{6}$  & $  1.7 \times 10^{7}$  \\
\end{tabular}
\vspace{0.1cm}
\caption{Optimal number of samples per level for the multilevel \name for example 1 with 
control-variate technique.}
\label{Samples per level for the drift implicit MLMC tau-leap for the decaying-dimerizing example control_variate} 
\end{table}

\begin{table}[H]
\centering
\begin{tabular}{l | c | c | c | c }
 tol  & $0.5$ &  $0.2$ &  $0.1$ & $0.05 $   \\
\hline \hline
 Needed number of samples  & $ 1.3 \times 10^{5}$ & $ 10^{6}  $ & $ 4.2 \times 10^{6}   $  &    $1.7 \times 10^{7}$
  \\
\end{tabular}
\vspace{0.1cm}
\caption{Needed number of samples for the MC \name for Example 1.}
\label{Samples for the drift implicit MC tau-leap method for the decaying-dimerizing example}
\end{table}

\bibliographystyle{abbrv}

\end{document}